\documentclass[twoside, 10pt]{article}

\usepackage{amsmath, amscd, amsfonts, amssymb, amsthm, latexsym, url, color}

\usepackage{graphicx}

\usepackage[left=1.3in, right=1.2in, top=1in, bottom=0.8in, includefoot, headheight=13.6pt]{geometry}

\input{xy}
\xyoption{all}

\usepackage{fancyhdr}

\usepackage{sectsty}
\sectionfont{\Large\sc\centering}
\subsectionfont{\sc}

\newtheorem{theorem}{Theorem}[section]
\newtheorem{lemma}[theorem]{Lemma}
\newtheorem{corollary}[theorem]{Corollary}
\newtheorem{proposition}[theorem]{Proposition}

\theoremstyle{definition}
\newtheorem{definition}[theorem]{Definition}
\newtheorem{remark}[theorem]{Remark}
\newtheorem{example}[theorem]{Example}

\newcommand{\al}{\alpha}
\newcommand{\bb}{\mathbb}

\newcommand{\bor}{\partial}

\newcommand{\Cl}[1]{\overline{\{#1\}}}
\newcommand{\comment}[1]{}

\newcommand{\dblecurly}[1]{\{\!\{ #1 \}\!\}}

\newcommand{\into}{\hookrightarrow}
\newcommand{\isoto}{\stackrel{\simeq}{\to}}

\newcommand{\mult}[1]{#1^{\!\times}}

\newcommand{\omegaa}{\mbox{\begin{large}$\omega$\end{large}}}
\newcommand{\onto}{\twoheadrightarrow}
\newcommand{\op}{\operatorname}
\newcommand{\res}{\overline}
\newcommand{\roi}{\mathcal{O}}

\newcommand{\sub}[1]{{\mbox{\scriptsize #1}}}

\newcommand{\To}{\longrightarrow}
\newcommand{\ul}[1]{\underline{#1}}

\newcommand{\xto}{\xrightarrow}

\renewcommand{\cal}{\mathcal}
\renewcommand{\hat}{\widehat}
\renewcommand{\frak}{\mathfrak}
\newcommand{\indlim}{\varinjlim}
\renewcommand{\tilde}{\widetilde}

\renewcommand{\ker}{\operatorname{Ker}}
\renewcommand{\projlim}{\varprojlim}

\DeclareMathOperator{\Frac}{Frac}

\DeclareMathOperator{\RES}{res}

\DeclareMathOperator*{\rprod}{\prod\nolimits^{\prime}\hspace{-1mm}}
\DeclareMathOperator{\RRES}{Res}

\DeclareMathOperator{\Spec}{Spec}

\DeclareMathOperator{\Tr}{Tr}

\renewcommand{\omegaa}{\mbox{\begin{Large}$\omega$\end{Large}}}

\begin{document}
\title{Grothendieck's trace map for arithmetic surfaces via residues and higher adeles}

\author{\sc Matthew Morrow}

\date{}



\maketitle

\begin{abstract}
We establish the reciprocity law along a vertical curve for residues of differential forms on arithmetic surfaces, and describe Grothendieck's trace map of the surface as a sum of residues. Points at infinity are then incorporated into the theory and the reciprocity law is extended to all curves on the surface. Applications to adelic duality for the arithmetic surface are discussed.


\end{abstract}

\section{Introduction}
Grothendieck's trace map for a smooth, projective curve over a finite field can be expressed as a sum of residues over all closed points of the curve (see \cite[III.7.14]{Hartshorne1977}). This result was generalized to algebraic surfaces by A.~Parshin \cite{Parshin1976} using his theory of two-dimensional adeles and residues for two-dimensional local fields. The theory for arbitrary dimensional algebraic varieties is essentially contained in A.~Beilinson's short paper on higher dimensional adeles \cite{Beilinson1980}, with considerable additional work by J.~Lipman \cite{Lipman1984}, V.~Lomadze \cite{Lomadze1981}, D.~Osipov \cite{Osipov1997}, A.~Yekutieli \cite{Yekutieli1992}, et al. In all these existing cases one restricts to varieties over a field. The purpose of this and the author's earlier paper \cite{Morrow2009} is to provide the first extension of the theory to non-varieties, namely to arithmetic surfaces, even taking into account the points `at infinity'. 

In the standard approach to Grothendieck duality of algebraic varieties using residues, there are three key steps. Firstly one must define suitable local residue maps, either on spaces of differential forms or on local cohomology groups (the latter approach is followed by E.~Kunz \cite{Kunz2008} using Grothendieck's residue symbol \cite[III.\S9]{Hartshorne1966}). Secondly, the local residue maps are used to define the dualizing sheaf, and finally the local residue maps must be patched together to define Grothendieck's trace map on the cohomology of the dualizing sheaf. In a previous paper, the author \cite{Morrow2009} carried out most of the first two steps for arithmetic surfaces, as we now explain.

Section \ref{section_review} provides a detailed summary of the required results from \cite{Morrow2009}, while also establishing several continuity and vanishing results which are required later. Briefly, given a two-dimensional local field $F$ of characteristic zero and a fixed local field $K\le F$, we introduced (see section \ref{ss_two_dim_local_fields}) a relative residue map \[\RRES_F:\Omega_{F/K}^\sub{cts}\to K,\] where $\Omega_{F/K}^\sub{cts}$ is a suitable space of `continuous' relative differential forms. In the case $F\cong K((t))$, this is the usual residue map; but if $F$ is of mixed characteristic, then this residue map is new (though versions of it appear in I.~Fesenko's two-dimensional adelic analysis \cite[\S27, Prop.]{Fesenko2008a} and in D.~Osipov's geometric counterpart \cite[Def.~5]{Osipov1997} to this paper). Then the reciprocity law for two-dimensional local rings was proved, justifying our definition of the relative residue map for mixed characteristic fields. For example, suppose $A$ is a characteristic zero, two-dimensional, normal, complete local ring with finite residue field, and fix the ring of integers of a local field $\roi_K\le A$. To each height one prime $y\subset A$ one associates the two-dimensional local field $\Frac\hat{A_y}$ and thus obtains a residue map $\RRES_y:\Omega_{\Frac A/K}^1\to K$ (see section \ref{ss_two_dim_complete_local_ring}). We showed \[\sum_y\RRES_y\omega=0\] for all $\omega\in\Omega_{\Frac A/K}^1$. The main new result in section \ref{section_review} is lemma \ref{lemma_continuous_wrt_m_adic_topology}, stating that the residue map $\RRES_y$ is continuous with respect to the $\frak{m}$-adic topology on $A$.

Geometrically, if $\pi:X\to\Spec\roi_K$ is an arithmetic surface and one chooses a closed point $x\in X$ and an irreducible curve $y\subset X$ passing through $x$, then one obtains a residue map \[\RRES_{x,y}:\Omega_{K(X)/K}^1\to K_{\pi(x)},\] where $K_{\pi(x)}$ is the completion of $K$ at the prime sitting under $x$ (see section \ref{subsection_geometric_case} for details). The established reciprocity law now takes the following form: \[\sum_{y\sub{ s.t. }y\ni x}\RRES_{x,y}\omega=0,\] where one fixes $\omega\in\Omega_{K(X)/K}^1$ and the summation is taken over all curves $y$ passing through a fixed point $x$.

As discussed above, the second step in a residue-theoretic approach to Grothendieck duality is a suitable description of the dualizing sheaf. This was also given in \cite{Morrow2009}: if $\pi:X\to\Spec\roi_K$ is an arithmetic surface (the precise requirements are those given at the start of section \ref{section_reciprocity_along_vertical_curves}), then the dualizing sheaf $\omegaa_\pi$ of $\pi$ can be described as follows: \begin{align*}\omegaa_\pi(U)=\{\omega\in\Omega_{K(X)/K}^1:&\RRES_{x,y}(f\omega)\in\hat{\roi_{K,\pi(x)}}\mbox{ for}\\ &\mbox{all } x\in y\subset U\mbox{ and }f\in\roi_{X,y}\}\end{align*} where $x$ runs over all closed points of $X$ inside $U$ and $y$ runs over all curves containing $x$. 

This paper treats the third step of the process. In order to patch the local residues together to define the trace map on cohomology, one must, just as in the basic case of a smooth, projective curve, establish certain reciprocity laws. For an arithmetic surface, these take the form: \[\sum_{y\sub{ s.t. }y\ni x}\RRES_{x,y}\omega=0,\quad\quad\sum_{x\sub{ s.t. }x\in y}\RRES_{x,y}\omega=0.\] In both cases one fixes $\omega\in\Omega_{K(X)/K}^1$, but the first summation is taken over all curves passing through a fixed point $x$ while the second summation is over all closed points of a fixed vertical curve $y$. The first of these laws, namely reciprocity around a point, has already been discussed, while section \ref{section_reciprocity_along_vertical_curves} establishes the reciprocity law along a vertical curve: the key idea of the proof is to reduce to the case when $\roi_K$ is a complete discrete valuation ring and then combine the reciprocity law around a point with the usual reciprocity law along the generic fibre.

Section \ref{section_trace_map_via_residues} uses the Parshin-Beilinson higher adeles for coherent sheaves to express Grothendieck's trace map \[\mbox{tr}_\pi:H^1(X,\omegaa_\pi)\to\roi_K\] as a sum of the residue maps $(\RRES_{x,y})_{x,y}$. Indeed, the reciprocity laws imply that our residue maps descend to cohomology: the argument is analogous to the case of a smooth, projective curve, except we must work with adeles for two-dimensional schemes rather than the more familiar adeles of a curve. Remark \ref{remark_higher_dimension} explains the basic framework of the theory in arbitrary dimensions.

Whereas the material discussed above is entirely scheme-theoretic, the final part of the paper is the most important and interesting from an arithmetic perspective as it incorporates archimedean points (points at infinity). It is natural to ask whether there exists a reciprocity law for {\em all} curves on $X$, not merely the vertical ones, when $\roi_K$ is the ring of integers of a number field. By compactifying $\Spec\roi_K$ and $X$ to include archimedean points in section \ref{section_archimedean}, we indeed prove a reciprocity law for any horizonal curve $y$ on $X$. Owing to the non-existence (at least naively) of $\Spec\bb{F}_1$, this takes the form \[\prod_{x\sub{ s.t. }x\in y}\psi_{x,y}(\omega)=1,\] where $\psi_{x,y}:\Omega_{K(X)/K}^1\to S^1$ are {\em absolute residue maps}, i.e.~additive characters, and $\omega\in\Omega_{K(X)/K}^1$. This provides detailed proofs of various claims made in \cite[\S27, \S28]{Fesenko2008a} concerning the foundations of harmonic analysis and adelic duality for arithmetic surfaces, and extends Parshin's absolute reciprocity laws for algebraic surfaces to the arithmetic case. Essentially this yields a framework which encodes both arithmetic duality of $K$ and Grothendieck duality of $X\to S$, and which would be equivalent to Serre duality were $X$ a geometric surface; a comparison of these results with Arakelov theory has yet to be carried out but there is likely an interesting connection.

Combined with \cite{Morrow2009}, which should be seen as a companion to this article and which contains a much more extensive introduction to the subject, these results provide a theory of residues and explicit duality for arithmetic surfaces. The analogous theory for an algebraic surface fibred smoothly over a curve is due to Osipov \cite{Osipov1997}, who proved, using Parshin's reciprocity laws for an algebraic surface, the analogues of our reciprocity laws around a point and along a vertical curve, and also showed that the sum of residues induces the trace map on cohomology.

\subsection{Notation}
When differential forms appear in this paper, they will be $1$-forms; so we write $\Omega_{A/R}$ in place of $\Omega_{A/R}^1$ to ease notation. $\Frac$ denotes the total ring of fractions; that is, if $R$ is a commutative ring then $\Frac R=S^{-1}R$, where $S$ is the set of regular elements in $R$. The maximal ideal of a local ring $A$ is usually denoted $\frak{m}_A$; an exception to this rule is when $A=\cal{O}_F$ is a discrete valuation ring with fraction field $F$, in which case we prefer the notation $\frak{p}_F$.

When $X$ is a scheme and $n\ge 0$, we write $X^n$ for the set of codimension $n$ points of $X$. $X_0$ denotes the closed points of $X$. $X$ will typically be two-dimensional, in which case we will often identify any $y\in X^1$ with the corresponding irreducible subscheme $\res{\{y\}}$; moreover, `$x\in y$' then more precisely means that $x$ is a codimension $1$ point of $\res{\{y\}}$. `Curve' usually means `irreducible curve'. Given $z\in X$, the maximal ideal of the local ring $\roi_{X,z}$ is written $\frak{m}_{X,z}$.

$I\subset^1\! A$ means that $I$ is a height one ideal of the ring $A$.

\subsection{Acknowledgements}
I am thankful to I.~Fesenko, A.~Beilinson, and A.~Yekutieli for discussions about this work. Parts of this research were funded by the Simons foundation and the EPSRC, and I am grateful to both organisations for their support.

I heartily thank the referee for suggestions and for reading the article in considerable depth, bringing to my attention several important issues which would have otherwise escaped my notice.

\section{Relative residue maps in dimension two}\label{section_review}

In \cite{Morrow2009}, a theory of residues on arithmetic surfaces was developed; we repeat here the main definitions and properties, also verifying several new results which will be required later.

\subsection{Two-dimensional local fields}\label{ss_two_dim_local_fields}
Suppose first that $F$ is a two-dimensional local field (i.e., a complete discrete valuation field whose residue field $\res{F}$ is a local field\footnote{In this paper our local fields always have finite residue fields, though many of the calculations continue to hold in the case of perfect residue fields.}) of characteristic zero, and that $K\le F$ is a local field (this local field $K$ will appear naturally in the geometric applications); write $\Omega_{F/K}^\sub{cts}=\Omega_{\roi_F/\roi_K}^\sub{sep}\otimes_{\roi_F}F$ (for a module over a local ring $A$, we write $M^\sub{sep}=M/\bigcap_{n\ge 0}\frak{m}_A^nM$ for the maximal separated quotient of $M$). Let $k_F$ be the algebraic closure of $K$ inside $F$; this is a finite extension of $K$ and hence is also a local field.

If $F$ has equal characteristic then any choice of a uniformiser $t\in F$ induces a unique $k_F$-isomorphism $F\cong k_F((t))$, and $\Omega_{\roi_F/\roi_K}^\sub{sep}=\roi_F\,dt$. The {\em relative residue map}, which does not depend on $t$, is the usual residue map which appears in the theory of curves over a field (e.g. \cite[II.7]{Serre1988}): \[\RES_F:\Omega_{F/K}^\sub{cts}\to k_F,\quad f\,dt\mapsto\mbox{coeft}_{t^{-1}}f,\] where the notation means that $f$ is to be expanded as a series in powers of $t$ and the coefficient of $t^{-1}$ is to be taken.

If $F$ is a mixed characteristic two-dimensional local field then $F/k_F$ is an infinite extension of complete discrete valuation fields, and $F$ is called {\em standard} if and only if $e(F/k_F)=1$. If $F$ is standard then any choice of a first local parameter $t\in\roi_F$ (i.e., $\res{t}$ is a uniformiser in the local field $\res{F}$) induces a unique $k_F$-isomorphism $F\cong k_F\dblecurly{t}$($:=$the completion of $\Frac(\roi_{k_F}[[t]])$ at the discrete valuation corresponding to the prime ideal $\frak{p}_{k_F}\roi_{k_F}[[t]]$; see \cite[Ex.~2.10]{Morrow2009}), and $\Omega_{\roi_F/\roi_K}^\sub{sep}=\roi_F\,dt$; so we may define \[\RES_F:\Omega_{F/K}^\sub{cts}\to k_F,\quad f\,dt\mapsto-\mbox{coeft}_{t^{-1}}f,\] which was shown in \cite[Prop.~2.19]{Morrow2009} not to depend on the choice of $t$. (The notation again means that $f$ is to be expanded as a series in powers in $t$, but this time in the field $k_F\dblecurly{t}$, and the coefficient of $t^{-1}$ taken). If $F$ is not necessarily standard, then choose a subfield $M\le F$ which is a standard two-dimensional local field, such that $F/M$ is a finite extension, and which satisfies $k_M=k_F$. The {\em relative residue map} in this case is defined by \[\RES_F=\RES_M\circ\Tr_{F/M}:\Omega_{F/K}^\sub{cts}\to k_F,\] which was shown in \cite[Lem.~2.21]{Morrow2009} not to depend on $M$.

In both cases, it is also convenient to write $\RRES_F=\Tr_{k_F/K}\circ\RES_F:\Omega_{F/K}^\sub{cts}\to K$. Also note that $\RES_F$ is $k_F$-linear, and that therefore $\RRES_F$ is $K$-linear. The expected functoriality result holds:

\begin{lemma}
Let $L$ be a finite extension of $K$. Then $\Omega_{L/K}^\sub{cts}$ is naturally isomorphic to $\Omega_{F/K}^\sub{cts}\otimes_F L$, so that there is a trace map $\Tr_{L/F}:\Omega_{L/K}^\sub{cts}\to\Omega_{F/K}^\sub{cts}$. If $\omega\in\Omega_{L/K}^\sub{cts}$, then \[\RRES_F(\Tr_{L/F}\omega)=\RRES_L\omega\] in $K$.
\end{lemma}
\begin{proof}
In the equal characteristic case this is classical; e.g., see \cite[II.12 Lem.~5]{Serre1988}. For the mixed characteristic case, see \cite[Prop.~2.22]{Morrow2009}.
\end{proof}

Next we show a couple of results on the continuity of residues which, though straightforward, will be frequently employed. A stronger, similar result is lemma \ref{lemma_continuous_wrt_m_adic_topology} later.

\begin{lemma}\label{lemma_residues_of_integral_forms}
Suppose that $\omega\in\Omega_{F/K}^\sub{cts}$ is integral; i.e., belongs to the image of $\Omega_{\roi_F/\roi_K}^\sub{sep}$. Then $\RES_F\omega\in\roi_{k_F}$ and so $\RRES_F\omega\in\roi_K$; in fact, if $F$ is equal characteristic, then $\RES_F\omega=0$.
\end{lemma}
\begin{proof}
In the equal characteristic or standard case this follows immediately from the definitions. In the non-standard, mixed characteristic case, one picks a standard subfield $M$ as above and uses a classical formula for the different of $\roi_F/\roi_M$ to show that the trace map $\Omega_{F/K}^\sub{cts}\to\Omega_{M/K}^\sub{cts}$ may be pulled back to $\Omega_{\roi_F/\roi_K}^\sub{sep}\to\Omega_{\roi_M/\roi_K}^\sub{sep}$, from which the result follows. See \cite[\S2.3.4]{Morrow2009} for the details.
\end{proof}

\begin{remark}\label{remark_commutativity_of_residues_with_reduction}
It was also shown in \cite[Cor.~2.23]{Morrow2009} that, when $F$ has mixed characteristic, the following diagram commutes:
\[\begin{CD}
\Omega_{\roi_F/\roi_K}^\sub{sep}& @>\RRES_F>> &\roi_{K}\\
@VVV&&@VVV\\
\Omega_{\res{F}/\res{K}}&@>>e(F/K)\RRES_{\res{F}}>&\res{K}\\
\end{CD}\]
The top horizontal arrow here makes sense by the previous lemma, and the lower horizontal arrow is the ramification degree $e(F/K)$ times the residue map for the local field $\res{F}$ of finite characteristic, which contains the finite field $\res{K}$.
\end{remark}

\begin{corollary}\label{corollary_continuity_of_residues_wrt_valuation_topology}
Fix $\omega\in\Omega_{F/K}^\sub{cts}$. Then \[F\to K,\quad f\mapsto \RRES_F(f\omega)\] is continuous with respect to the discrete valuation topologies on $F$ and $K$; in fact, if $F$ is equal characteristic, then it is even continuous with respect to the discrete topology on $K$.
\end{corollary}
\begin{proof}
After multiplying $\omega$ by a non-zero element of $F$, we may assume that $\omega$ is integral in the sense of the previous lemma. If $F$ is equal characteristic then $\ker(f\mapsto \RES_F(f\omega))$ contains the open set $\roi_F$, proving continuity with respect to the discrete topology on $K$. Now assume $F$ has mixed characteristic and let $\pi$ be a uniformiser of $K$; since $F/K$ is an extension of complete discrete valuation fields, we may put $e=e(F/K)=\nu_F(\pi)>0$. Then the previous lemma implies \[\RRES(\frak{p}_F^{es}\omega)=\RRES(\pi^s\roi_F\omega)=\pi^s\RRES(\roi_F\omega)\subseteq\frak{p}_K^s\] for all $s\in\bb{Z}$, proving continuity with respect to the discrete valuation topologies.
\end{proof}

\subsection{Two-dimensional complete rings}\label{ss_two_dim_complete_local_ring}
Let $A$ be a two-dimensional, normal, complete, local ring of characteristic zero, with a finite residue field of characteristic $p$; set $F=\Frac A$. Then there is a unique ring homomorphism $\bb{Z}_p\to A$ and it is a closed embedding; let $\roi_K$ be a finite extension of $\bb{Z}_p$ inside $A$; i.e., $\roi_K$ is the ring of integers of $K$, which is a finite extension of $\bb{Q}_p$.

If $y\subset A$ is a height one prime (we often write $y\subset^1\!A$), then $\hat{A_y}$ is a complete discrete valuation ring; its field of fractions $F_y:=\Frac\hat{A_y}$ is a two-dimensional local field containing $K$. Moreover, there is a natural isomorphism $\Omega_{A/\roi_K}^\sub{sep}\otimes_A \hat{A_y}\cong \Omega_{\hat{A_y}/K}^\sub{sep}$ (see \cite[Lem.~3.8]{Morrow2009}); so we define $\RRES_y:\Omega_{A/\roi_K}^\sub{sep}\otimes_A F\to K$ to be the composition \[\Omega_{A/\roi_K}^\sub{sep}\otimes_A F\To \Omega_{A/\roi_K}^\sub{sep}\otimes_A F_y\cong\Omega_{F_y/K}^\sub{cts}\xto{\RRES_{F_y}} K.\] The definition of the residue maps is justified by the following reciprocity law:

\begin{theorem}\label{theorem_reciprocity_for_2d_complete_rings}
Let $\omega\in\Omega_{A/\roi_K}^\sub{sep}\otimes_A F$; then for all but finitely many height one primes $y\subset A$ the residue $\RRES_y\omega$ is zero, and \[\sum_{y\subset^1\!A}\RRES_y\omega=0.\]
\end{theorem}
\begin{proof}
See \cite[Thm.~3.10]{Morrow2009}.
\end{proof}

As is often the case, the residue law was reduced to a special case by taking advantage of functoriality:

\begin{lemma}\label{lemma_functoriality_for_2d_complete_rings}
Suppose that $C$ is a finite extension of $A$ which is also normal; set $L=\Frac C$. Then for any $\omega\in\Omega_{C/\roi_K}^\sub{sep}\otimes_C L$ and any height one prime $y\subset A$, we have \[\RRES_y(\Tr_{L/F}\omega)=\sum_{Y|y}\RRES_Y\omega,\] where $Y$ varies over the finitely many height one primes of $C$ which sit over $y$.  
\end{lemma}
\begin{proof}
See \cite[Thm.~3.9]{Morrow2009}.
\end{proof}

The proof of the reciprocity theorem also required certain results on the continuity of the residues whose proofs were omitted in \cite{Morrow2009}; we shall require similar such results several times in this article and now is a convenient opportunity to establish them:

\begin{lemma}
Set $B=\roi_K[[t]]$, $M=\Frac{B}$ and let $\omega\in\Omega_{B/\roi_K}^\sub{sep}\otimes_B M$; then, for any height one prime $y\subset B$, the map \[B\to K,\quad f\mapsto \RRES_yf\omega\] is continuous with respect to the $\frak{m}_B$-adic topology on $B$ and the discrete valuation topology on $K$.
\end{lemma}
\begin{proof}
We first consider the case when $y=\rho B$ is generated by an irreducible Weierstrass polynomial $\rho(t)\in\roi_K[t]$. Let $K'$ be a sufficiently large finite extension of $K$ such that $\rho$ splits into linear factors in $K'$; the decomposition has the form $\rho(t)=\prod_{i=1}^d(t-\lambda_i)$ with $d=\deg \rho$ and $\lambda_i\in\frak{p}_{K'}$ since $h$ is a Weierstrass polynomial. Put $B'=\roi_{K'}[[t]]$ and $M'=\Frac B'$. According to functoriality of residues (the previous lemma), we have \[\RRES_y \Tr_{M'/M}\omega=\sum_{i=1}^d\RRES_{Y_i}\omega\] for all $\omega\in\Omega_{B'/\roi_K}^\sub{sep}\otimes_{B'} M'$, where $Y_i=(t-\lambda_i)B'$. Since multiplication by $f\in B$ commutes with the trace map, it is now enough to prove that \[B'\to K,\quad f\mapsto \RRES_{Y_i}f\omega\] is continuous for all $i$ and all $\omega\in\Omega_{B'/\roi_K}^\sub{sep}\otimes_{B'} M'$. In other words, replacing $K$ by $K'$ and $B$ by $B'$, we have reduced to the case when $\rho(t)$ is a linear polynomial: $\rho(t)=t-\lambda$, with $\lambda\in\frak{p}_K$. After another reduction, we will prove the continuity claim in this case.

Let $\pi$ be a uniformiser for $K$. It is well-known that $\Omega_{B/\roi_K}^\sub{sep}=B\,dt$ and that any element of $M$ can be written as a finite sum of terms of the form \[\frac{\pi^n g}{h^r},\] with $h\in\roi_K[t]$ an irreducible Weierstrass polynomial, $r>0$, $n\in\bb{Z}$, and $g\in B$ (a proof was given in \cite[Lem.~3.4]{Morrow2009}). By continuity of addition $K\times K\xto{+}K$ and of the multiplication maps $B\xto{\times g}B$, $K\xto{\times\pi^n}K$, it is enough to treat the case \[\omega=h^{-r}\,dt,\] where $h\in\roi_K[t]$ is an irreducible Weierstrass polynomial.

Now return to $y=\rho B$, $\rho=t-\lambda$. If $h\neq\rho$, then $h^{-r}dt\in\Omega_{B/\roi_K}^\sub{sep}\otimes_B B_y$, and so $\RRES_y(B\omega)=0$ by lemma \ref{lemma_residues_of_integral_forms}, which is certainly enough. Else $h=\rho$, which we now consider. To obtain more suggestive notation, we write $t_y:=\rho(t)=t-\lambda$; thus \[\omega=h^{-r}\,dt=t_y^{-r}\,dt_y.\] Let $m\ge 0$; we claim that if $n\ge m+r$ then $\RRES_y(\frak{m}_B^n\omega)\subseteq\frak{p}_K^m$. Since $\lambda$ is divisible by $\pi$, the maximal ideal of $B$ is generated by $\pi$ and $t_y$: \[\frak{m}_B=\langle \pi,t\rangle=\langle \pi,t_y\rangle.\] Therefore an arbitrary element of $\frak{m}_B^n$ is a sum of terms of the form $\pi^\al t_y^\beta g$, with $g\in B$, $\al,\beta\ge 0$, and $\al+\beta\ge n$, and so it is enough to consider such an element. Moreover, again since $\pi$ divides $\lambda$, there is a unique continuous isomorphism \[\roi_K[[t_y]]\isoto\roi_K[[t]],\quad t_y\mapsto t-\lambda,\] and therefore $g\in B$ may be written as \[g=\sum_{j=0}^{r-1}a_jt_y^j+t_y^r g_1\] with $a_j\in \roi_K$ and $g_1\in B$ (we could extend this expansion to infinity, of course, but since we are trying to prove continuity, it is better not to risk confusion between `formal series' and `convergent series'). Then \[\RRES_y (\pi^\al t_y^\beta g\omega)=\pi^{\al}\RRES_y\left(t_y^{\beta-r}\sum_{j=0}^{r-1}a_jt_y^j\,dt_y\right) + \pi^\al \RRES_y(t_y^\beta g_1\,dt_y).\tag{\dag}\] The second residue is zero by lemma \ref{lemma_residues_of_integral_forms} again since $t_y^{\beta}g_1\in B$. If $\beta\ge r$ then the first residue is zero for the same reason; but if $\beta <r$ then it follows that $\al>m$, whence the first residue is $\pi^\al a_{r-\beta-1}\in\frak{p}_K^\al\subseteq\frak{p}_K^m$. So in any case, (\dag) belongs to $\frak{p}_K^m$, completing the proof of our claim and thereby showing the desired continuity result for $y=\rho B$.

Having treated the case of a prime $y$ generated by a Weierstrass polynomial, we must secondly consider $y=\pi B$. By exactly the same argument as above, we may assume that $\omega=h^{-r}\,dt$, with $h$ an irreducible Weierstrass polynomial. Then $M_y=K\dblecurly{t}$ and $h^{-r}\in B_y$; hence $h^{-r}$ may be written as a series \[h^{-r}=\sum_{j\in\bb{Z}} a_jt^j\in\roi_K\dblecurly{t}\] where $a_j\to0$ in $\roi_K$ as $j\to-\infty$. Let $m\ge 0$ be fixed, and pick $J>2$ such that $a_j\in\frak{p}_K^m$ whenever $j\le -J$. We claim that if $n\ge J-2+m$ then $\RRES_y(\frak{m}_B^n\omega)\subseteq\frak{p}_K^m$. Since an arbitrary element of $\frak{m}_B^n$ is a sum of terms of the form $\pi^\al t^\beta g$, with $g\in B$, $\al,\beta\ge 0$, and $\al+\beta\ge n$, it is enough it consider such an element; write $g=\sum_{i=0}^\infty b_it^i$. Then
\begin{align*}
\RRES_y (\pi^\al t^\beta g\omega)
	&=\RRES_y (\pi^\al t^\beta gh^{-r}\,dt)\\
	&=-\pi^\al\mbox{coeft}_{t^{-1}}\left(t^\beta\sum_{i=0}^\infty b_it^i\sum_{j\in\bb{Z}} a_jt^j\right)\\
	&=-\pi^\al\sum_{i=0}^\infty b_i a_{-i-\beta-1}\\
	&\in\begin{cases}\frak{p}_K^{\al+m} & \mbox{if }\beta\ge J-2 \\ \frak{p}_K^\al &\mbox{in any case}.\end{cases}
\end{align*}
But $\al+\beta\ge J-2+m$ and so if it is not the case that $\beta\ge J-2$, then it follows that $\al\ge m$; so, regardless of which inequality holds, we obtain $\RRES_y (\pi^\al t^\beta g\omega)\in\frak{p}_K^m$, as required.
\end{proof}

Now we extend the lemma to the general case of our two-dimensional, normal, complete, local ring $A$. This result is a significant strengthening of corollary \ref{corollary_continuity_of_residues_wrt_valuation_topology}, since the $\frak{m}_A$-adic topology on $A$ is considerably finer than the $y$-adic topology, for any $y\subset^1\! A$.

\begin{lemma}\label{lemma_continuous_wrt_m_adic_topology}
Let $\omega\in\Omega_{A/\roi_K}^\sub{sep}\otimes_A F$; then, uniformly in $y$, the map \[A\to K,\quad f\mapsto \RRES_yf\omega\] is continuous with respect to the $\frak{m}_A$-adic topology on $A$ and the discrete valuation topology on $K$.
\end{lemma}
\begin{proof}
Firstly, it is enough to prove that the given map is continuous for any fixed $y$; the uniformity result then follows from the fact that, for almost all $y\subset^1\!A$, $\omega$ belongs to $\Omega_{A_y/\roi_K}^\sub{sep}$ and $y$ does not contain $\frak{p}_K$; for such primes, $\RRES_yA\omega=0$ by lemma \ref{lemma_residues_of_integral_forms}.

By Cohen structure theory \cite{Cohen1946} (the details of the argument are in \cite[Lem.~3.3]{Morrow2009}), there is a subring $B\le A$ containing $\roi_K$ which is isomorphic to $\roi_K[[t]]$ and such that $A$ is a finitely-generated $B$-module; set $M=\Frac B$. Write $\omega=g\omega_0$ for some $g\in F$ and $\omega_0\in\Omega_{B/\roi_K}^\sub{sep}\otimes_B M$.

Now we make some remarks on continuity of the trace map. $\Tr_{F/M}(Ag)$ is a finitely generated $B$-module and so there exists $g_0\in\mult{M}$ such that $\Tr_{F/M}(Ag)\subseteq Bg_0$. Moreover, since $A/B$ is a finite extension of local rings, one has $\frak{m}_A^s\subseteq\frak{m}_BA$ for some $s>0$. Hence $\Tr_{F/M}(\frak{m}_A^{ns}g)\subseteq\frak{m}_B^ng_0$ for all $n\ge 0$, meaning that the restriction of the trace map to $Ag\to Bg_0$ is continuous with respect to the $\frak{m}$-adic topologies on each side. It immediately follows that \[\tau:A\to B,\quad f\mapsto\Tr_{F/M}(fg)g_0^{-1}\] is both well-defined and continuous.

Functoriality (lemma \ref{lemma_functoriality_for_2d_complete_rings}) implies that for any $y\subset^1\! B$, \[\sum_{Y|y}\RRES_Y f\omega=\RRES_y\Tr_{F/M}(f\omega)\] for all $f\in A$, where $Y$ varies over the finitely many height one primes of $A$ which sit over $y$. The right hand side may be rewritten as \[\RRES_y(\tau(f)\,g_0\omega_0)\] where $g_0\omega_0\in\Omega_{B/\roi_K}^\sub{sep}\otimes_B M$; according to the previous lemma, this is a continuous function of $f$. In conclusion, \[A\to K,\quad f\mapsto\sum_{Y|y}\RRES_Y f\omega\tag{\dag}\] is continuous, which we will now use to show that each map $f\mapsto \RRES_Yf\omega$ is individually continuous, thereby completing the proof. Fix $m\ge 0$.

Let $Y_1,\dots,Y_l$ be the height one primes of $A$ sitting over $y$, and let $\nu_1,\dots,\nu_l$ denote the corresponding discrete valuations of $F$. If $l=1$ then there is nothing more to show, so assume $l>1$. Since the map \[F_{Y_i}\to K,\quad f\mapsto\RRES_{Y_i}(f\omega)\] is continuous with respect to the discrete valuation topologies on each side (corollary \ref{corollary_continuity_of_residues_wrt_valuation_topology}), there exists $S>0$ (which we may obviously assume is independent of $i$) such that $\RRES_{Y_i}(f\omega)\subseteq\frak{p}_K^m$ whenever $\nu_i(f)\ge S$. According to the approximation theorem for discrete valuations, there exists an element $e\in F$ which satisfies $\nu_1(e-1)\ge S$ and $\nu_i(e)\ge S$ for $i=2,\dots, l$. Now, since (\dag) remains continuous if we replace $\omega$ by $e\omega$, there also exists $J>0$ such that $\sum_{Y|y}\RRES_Y (fe\omega)\in\frak{p}_K^m$ whenever $f\in\frak{m}_A^J$.

So, if $f\in\frak{m}_A^J$ then  \[\RRES_{Y_1}(f\omega)=\RRES_{Y_1}(f(1-e)\omega)-\sum_{i=2}^l\RRES_{Y_i}(fe\omega)+\sum_{i=1}^l\RRES_{Y_i}(fe\omega)\] belongs to $\frak{p}_K^m$ since $\nu_1(f(1-e))\ge S$ and $\nu_i(fe)\ge S$ for $i=2,\dots,l$. That is, $\RRES_{Y_1}(\frak{m}_A^J\omega)\subseteq\frak{p}_K^m$, which proves the desired continuity result.
\end{proof}

\begin{remark}\label{remark_continuous_wrt_two_dim_topology}
The previous lemma can be reformulated as saying that the residue map $\RRES_{F_y}:\Omega_{F_y/K}^\sub{cts}\to K$ is continuous with respect to the valuation topology on $K$ and the vector space topology on $\Omega_{F_y/K}^\sub{cts}$, having equipped $F_y$ with its two-dimensional local field topology \cite{Madunts1995}.
\end{remark}

Finally, regarding vanishing of the residue of a differential form:

\begin{lemma}
Suppose that $\omega\in\Omega_{A/\roi_K}^\sub{sep}\otimes_A F$ is integral, in the sense that it belongs to the image of $\Omega_{A/\roi_K}^\sub{sep}$, and let $y\subset^1\! A$. Then $\RRES_y\omega\in\frak{p}_K$. If $y$ does not contain $p$ or if $y$ is the only height one prime of $A$ containing $p$, then $\RRES_y\omega=0$. 
\end{lemma}
\begin{proof}
If $y$ does not contain $p$ then $F_y$ is equal characteristic and we have already proved a stronger result in lemma \ref{lemma_residues_of_integral_forms}: $\RRES_y$ vanishes on the image of $\Omega_{A/\roi_K}^\sub{sep}\otimes_AA_y$. If instead $y$ is the only height one prime of $A$ containing $p$, then the vanishing claim follows from the reciprocity law and the previous case.

Finally, suppose $y$ contains $p$ but do not assume that it is the only height one prime to do so. Using functoriality of differential forms and remark \ref{remark_commutativity_of_residues_with_reduction}, we have a commutative diagram
\[\begin{CD}
\Omega_{A/\roi_K}^\sub{sep} & @>>> &\Omega_{\roi_{F_y}/\roi_K}^\sub{sep}& @>\RRES_{F_y}>> &\roi_{K}\\
@VVV&&@VVV&&@VVV\\
\Omega_{(A/y)/\res{K}} & @>>> & \Omega_{\res{F}_y/\res{K}}&@>>e(F_y/K)\RRES_{\res{F}_y}>&\res{K}\\
\end{CD}\]
The residue map $\RRES_{\res{F}_y}$ on the characteristic $p$ local field $\res{F}_y$ vanishes on integral differential forms; since $A/y$ belongs to the ring of integers of $\res{F}_y$, it follows immediately from the diagram that $\RRES_y\omega\in\frak{p}_K$.
\end{proof}

\begin{example}
This example will show that the previous lemma cannot be improved. We consider the `simplest' $A$ in which $p$ splits. Set $B=\bb{Z}_p[[T]]$, with field of fractions $M$, and let $A=B[\al]$ where $\al$ is a root of $f(X)=X^2-TX-p$, with field of fractions $F$. Since $f(X)$ does not have a root in $B/TB=\bb{Z}_p$, it does not have a root in $B$, and so $F/M$ is a degree two extension. Since $A$ is a finitely generated $B$-module, it is also a two-dimensional, complete local ring, and we leave it to the reader to check that $A$ is regular, hence normal.

In $A$, $p$ completely splits as $p=\al(T-\al)$, and therefore, setting $y=\al A$, the natural map \[\bb{Q}_p\dblecurly{T}=M_{pB}\to F_y\] is an isomorphism. Indeed, $f(X)$ splits in the residue field $B_{pB}/pB_{pB}=\bb{F}_p((T))$ into distinct factors and so Hensel's lemma implies that $f(X)$ splits in $\hat{B_{pB}}$; i.e., $\al\in\hat{B_{pB}}\subset M_{pB}$.

One readily checks that $\al\equiv-pT^{-1}\op{mod}p^2$ in $\hat{B_{pB}}=\hat{A_y}$, and therefore \[\RRES_y(\al\,dT)\equiv-p\mod{p^2}.\] In particular, $\RRES_y(\al\,dT)\neq 0$ even though $\al\,dT$ is integral. 
\end{example}

\subsection{Two-dimensional, finitely generated rings}\label{subsection_non_complete_case}
Next suppose that $\roi_K$ is a Dedekind domain of characteristic zero and with finite residue fields, and that $B$ is a two-dimensional, normal, local ring, which we assume is the localization of a two-dimensional, finitely-generated $\roi_K$-algebra. Set $A=\hat{B_{\frak{m}_B}}$ and $s=\frak{m}_B\cap\roi_K$. Then $A$ satisfies all the conditions introduced at the start of the previous subsection and contains $\roi_s:=\hat{\roi_{K,s}}$, which is the ring of integers of the local field $K_s:=\Frac\hat{\roi_{K,s}}$. Moreover, there is a natural identification $\Omega_{B/\roi_K}\otimes_B A=\Omega_{A/\roi_s}^\sub{sep}$ (see \cite[Lem.~3.11]{Morrow2009}). For each height one prime $y\subset B$, we may therefore define \[\RRES_y:\Omega_{\Frac B/K}\to K_s\] to be the composition \[\begin{CD}\Omega_{\Frac B/K}&@>>>& \Omega_{\Frac B/K}\otimes_{\Frac B} \Frac A\cong\Omega_{A/\roi_s}^\sub{sep}\otimes_A\Frac A&@>\sum_{y'|y}\RRES_{y'}>>& K_s\end{CD}\] where $y'$ varies over the finitely many primes of $A$, necessarily of height one, which sit over $y$.

The reciprocity law remains true in this setting:

\begin{theorem}
Let $\omega\in\Omega_{\Frac B/K}$; then for all but finitely many height one primes $y\subset B$ the residue $\RRES_y\omega$ is zero, and \[\sum_{y\subset^1\!B}\RRES_y\omega=0.\]
\end{theorem}
\begin{proof}
\cite[Thm.~3.13]{Morrow2009}
\end{proof}

The following vanishing identity will be useful:

\begin{lemma}\label{lemma_vanishing_of_residue}
Let $y\subset^1\! B$ and suppose that $\omega\in\Omega_{\Frac B/K}$ belongs to the image of $\Omega_{B_y/\roi_K}$. Then $\RRES_y\omega\in\roi_s$; in fact, if $y$ is horizontal (i.e., $y\cap\roi_K=0$) then $\RRES_y\omega=0$.

Secondly, suppose that there is only one height one prime $y$ of $B$ which is vertical (i.e., containing $s$) and that $\omega$ is in the image of $\Omega_{A/\roi_K}$. Then $\RRES_y\omega=0$.
\end{lemma}
\begin{proof}
The first claims follow from lemma \ref{lemma_residues_of_integral_forms}, since $y$ being horizontal is equivalent to the two-dimensional local fields $\Frac\hat{A_{y'}}$, with $y'\subset A$ sitting over $y$, being equi-characteristic. The second claim follows from the previous reciprocity law since any prime is either vertical or horizontal.
\end{proof}

\subsection{Geometrisation}\label{subsection_geometric_case}
Continue to let $\roi_K$ be a Dedekind domain of characteristic zero and with finite residue fields. Let $X$ be a two-dimensional, normal scheme, flat and of finite type over $S=\Spec\roi_K$, and let $\Omega_{X/S}=\Omega_{X/S}^1$ be the relative sheaf of one forms. Let $x\in X^2$ be a closed point sitting over a closed point $s\in S_0$, and let $y\subset X$ be an irreducible curve containing $x$. Identify $y$ with its local equation (i.e., corresponding prime ideal) $y\subset^1\!\roi_{X,x}$ and note that $\roi_{X,x}$ satisfies all the conditions which $B$ did in the previous subsection. Define the residue map $\RRES_{x,y}:\Omega_{K(X)/K}\to K_s\,(=\Frac\hat{\roi_{K,s}})$ to be \[\RRES_y:\Omega_{\Frac\roi_{X,x}/K}\To K_s.\] The reciprocity law now states that, for any fixed $\omega\in\Omega_{K(X)/K}$, \[\sum_{\substack{y\subset X\\\sub{s.t. }y\ni x}}\RRES_{x,y}\omega=0\] in $K_s$, where the sum is taken over all curves in $X$ which pass through $x$. For a few more details, see \cite[\S4]{Morrow2009}.

\section{Reciprocity along vertical curves}\label{section_reciprocity_along_vertical_curves}
As explained in the introduction, residues on a surface should satisfy {\em two} reciprocity laws, one as we vary curves through a fixed point, and another as we vary points along a fixed curve. The first was explained immediately above and now we will prove the second.

Let $\roi_K$ be a Dedekind domain of characteristic zero and with finite residue fields; denote by $K$ its field of fractions. Let $X$ be an $\roi_K$-curve; more precisely, $X$ is a normal scheme, proper and flat over $S=\Spec\roi_K$, whose generic fibre is a smooth, geometrically connected curve.

The aim of this section is to establish the following reciprocity law for vertical curves on an arithmetic surface:

\begin{theorem}\label{theorem_reciprocity_along_a_curve}
Let $\omega\in\Omega_{K(X)/K}$, and let $y\subset X$ be an irreducible component of a special fibre $X_s$, where $s\in S_0$. Then \[\sum_{x\in y}\RRES_{x,y}\omega=0\] in $K_s$, where the sum is taken over all closed points $x$ of $y$.
\end{theorem}

Here, as usual, $\roi_s=\hat{\roi_{K,s}}$ and $K_s=\Frac\roi_s$. The proof will consist of several steps. We begin with a short proof of a standard adelic condition:

\begin{lemma}\label{lemma_12_adelic_condition}
Let $y\subset X$ be an irreducible curve, let $f\in\roi_{X,y}$, and let $r\ge 1$. Then $f\in\roi_{X,x}+\frak{m}_{X,y}^r$ for all but finitely many closed points $x\in y$.

The result also holds after completion: if $f\in\hat{\roi_{X,y}}$, then $f\in\roi_{X,x}+\frak{m}_{X,y}^r\hat{\roi_{X,y}}$ for almost all $x$.
\end{lemma}
\begin{proof}
Let $U=\Spec A$ be an open affine neighbourhood of (the generic point of) $y$, let $\frak{p}\subset A$ be the prime ideal defining $y$, and set $P=A\cap\frak{p}^rA_\frak{p}$, $B=A/P$. If $b\in B$ is not a zero divisor, then $B/bB$ is zero-dimensional and so has only finitely many primes; hence only finitely many primes of $B$ contain $b$. Set \[\res{f}:=f\op{mod}\frak{m}_{X,y}^r\in A_\frak{p}/\frak{p}^rA_\frak{p}=\Frac{B};\] by what we have just proved, $\res{f}$ belongs to $B_{\frak{q}}$ for all but finitely many primes $\frak{q}\subset B$, i.e.~$f\in\roi_{X,x}+\frak{m}_{X,y}^r$ for all but finitely many $x\in y\cap U$. Since $U$ contains all but finitely many points of $y$, we have finished.

The complete version now follows from the identity $\hat{\roi_{X,y}}/\frak{m}_{X,y}^r\hat{\roi_{X,y}}=\roi_{X,y}/\frak{m}_{X,y}^r$.
\end{proof}

The lemma lets us prove that the theorem makes sense:

\begin{lemma}\label{lemma_convergence_of_residues_along_a_curve}
Let $\omega\in\Omega_{K(X)/K}$, and let $y\subset X$ be an irreducible component of a special fibre $X_s$, where $s\in S_0$. Then the sum $\sum_{x\in y}\RRES_{x,y}\omega$ converges in the valuation topology on $K_s$ (we will see that only countably many terms are non-zero).

Moreover, \[K(X)\to K_s,\quad h\mapsto\sum_{x\in y}\RRES_{x,y}(h\omega)\] is continuous with respect to the topology on $K(X)$ induced by the discrete valuation associated to $y$, and the $s$-adic topology on $K_s$.
\end{lemma}
\begin{proof}
For any point $z\in X$, let $\Omega_z$ denote the image of $\Omega_{\roi_{X,z}/\roi_K}$ inside $\Omega_{K(X)/K}$. Let $r\ge 0$.

Let $\pi\in\roi_K$ be a uniformiser at $s$, fix $\omega\in\Omega_{K(X)/K}$ and pick $a\ge 0$ such that $\pi^a\omega\in\Omega_y$. Then it easily follows from the previous lemma that, for any $r\ge 0$, $\pi^a\omega$ lies in $\Omega_x+\pi^r\Omega_y$ for almost all closed points $x\in y$. But lemma \ref{lemma_vanishing_of_residue} implies that if $x$ is any closed point of $y$ then $\RRES_{x,y}(\Omega_y)\subseteq\roi_s$, and moreover that if $x$ does not lie on any other irreducible component of the fibre $X_s$ then $\RRES_{x,y}(\Omega_x)=0$. We deduce that \[\RRES_{x,y}\pi^a\omega\in\pi^r\roi_s\] for almost all closed points $x\in y$. So $\RRES_{x,y}\omega\in\pi^{r-a}\roi_s$ for almost all $x\in y$; since this holds for all $r\ge 0$ we see that \[\sum_{x\in y}\RRES_{x,y}\omega\] converges and also that $\sum_{x\in y}\RRES_{x,y}\omega\in\pi^{-a}\roi_s$.

If $h\in K(X)$ satisfies $\nu(f)\ge b$ for some $b\in\bb{Z}$, then we may write $h=\pi^bu$ for some $u\in\roi_{X,y}$. This implies that $\pi^{a-b}h\omega\in\Omega_y$ and so, by what we have just shown, $\sum_{x\in y}\RRES_{x,y}h\omega\in\pi^{b-a}\roi_s$. This proves that $h\mapsto\sum_{x\in y}\RRES_{x,y}h\omega$ is continuous.
\end{proof}

\begin{remark}
The analogous vertical reciprocity law in the geometric setting is \cite[Prop.~6]{Osipov1997}, where Osipov gives an example to show that it really is possible for the sum of residues along the points of $y\subset X_s$ to contain infinitely many non-zero terms.
\end{remark}

We aim to reduce the vertical reciprocity law to the case of $\roi_K$ being a complete discrete valuation ring by using several lemmas on the functoriality of residues.

Let $s$ be a non-zero prime of $\roi_K$, and set $\roi_s=\hat{\roi_{K,s}}$, $K_s=\Frac{\roi_s}$ as usual. Set $\hat{X}=X\times_{\roi_K}\roi_s$ and let $p:\hat{X}\to X$ be the natural map. Then $p$ induces an isomorphism of the special fibres $\hat{X}_s\cong X_s$ and, for any point $x\in X_s$, $p$ induces an isomorphism of the completed local rings $\hat{\roi_{X,p(x)}}\cong \hat{\roi_{\hat{X},x}}$ (see e.g.~\cite[Lem.~8.3.49]{Liu2000}). From the excellence of $X$ it follows that $\roi_{\hat{X},x}$ is normal for all $x\in \hat{X}_s$, and therefore $\hat{X}$ is normal. So $\hat{X}$ is a $\roi_s$-curve, in the same sense as at the start of the section.

\begin{lemma}
Let $y\subset X$ be an irreducible curve and let $x$ be a closed point of $y$. Then the following diagram commutes:
\[\xymatrixcolsep{2.2cm}\xymatrix{
\Omega_{K(\hat{X})/K_s}\ar[dr]^{\sum_{y'|y}\RRES_{x',y'}}&\\
\Omega_{K(X)/K}\ar[r]_{\RRES_{x,y}}\ar[u]&\hat{K}\\
}\]
where $y'$ varies over the irreducible curves of $\hat{X}$ sitting over $y$ and $x'$ is the unique closed point sitting over $x$ (i.e.,~$p(x')=x$).
\end{lemma}
\begin{proof}
This essentially follows straight from the original definitions of the residue maps in subsections \ref{subsection_non_complete_case} and \ref{subsection_geometric_case}. Indeed, set $B=\roi_{X,x}$ and let $y\subset B$ be the local equation for $y$ at $x$, so that \[\RRES_{x,y}=\sum_{\substack{y''\subset\!^1\hat{B}\\y''|y}}\RRES_{y''}:\Omega_{\hat{B}/\roi_s}^\sub{sep} \otimes_{\hat{B}}\Frac\hat{B}\to K_s,\] where $y''$ varies over the height one primes of $\hat{B}$ sitting over $y$.

But we remarked above that there is a natural $\roi_s$-isomorphism $\hat{\roi_{\hat{X},x'}}\cong\hat{B}$, and this expression for the residues remains valid if $B$ is replaced by $\roi_{\hat{X},x'}$ and $y$ is replaced by some $y'$ sitting over $y$. Therefore
\begin{align*}
\RRES_{x,y}
	&=\sum_{\substack{y''\subset\!^1\hat{B}\\y''|y}}\RRES_{y''}\\
	&=\sum_{\substack{y'\subset\!^1\roi_{\hat{X},x'}\\y'|y}}\sum_{\substack{y''\subset\!^1\hat{B}\\y''|y'}}\RRES_{y''}\\
	&=\sum_{\substack{y'\subset\!^1\roi_{\hat{X},x'}\\y'|y}}\RRES_{y'}\\
	&=\sum_{y'|y}\RRES_{x',y'},
\end{align*}
as required.
\end{proof}

\begin{corollary}
Let $y\subset X$ be an irreducible component of the special fibre $X_s$ and let $x$ be a closed point of $y$; let $x'=p^{-1}(x)$, $y'=p^{-1}(y)$ be the corresponding point and curve on $\hat{X}_s\cong X_s$. Then the following diagram commutes:
\[\xymatrixcolsep{2.2cm}\xymatrix{
\Omega_{K(\hat{X})/K_s}\ar[dr]^{\RRES_{x',y'}}&\\
\Omega_{K(X)/K}\ar[r]_{\RRES_{x,y}}\ar[u]&K_s\\
}\]
Informally, this means that residues along the special fibre $X_s$ may be computed after completing $\roi_K$.
\end{corollary}
\begin{proof}
The unique irreducible curve of $\hat{X}$ sitting over $y$ is $y'$, so this follows from the previous lemma.
\end{proof}

\begin{corollary}
If the vertical reciprocity law holds for $\hat{X}/\roi_s$, then it holds for $X/\roi_K$.
\end{corollary}
\begin{proof}
This immediately follows from the previous corollary.
\end{proof}

In the remainder of the section (except remark \ref{remark_two_dim_vs_one_dim_residues}), we replace $X$ by $\hat{X}$ and $\roi_K$ by $\roi_s$, so that the base is a now a complete, discrete valuation ring (of characteristic zero, with finite residue field, with field of fractions $K$ being a local field).

The horizontal curves on $X$ are all of the form $\Cl{z}$ for a uniquely determined closed point $z$ of the generic fibre $X_\eta$. Moreover, because our base ring is now complete, $\Cl{z}$ meets the special fibre $X_s$ at a unique point $\frak{r}(z)$, which is necessarily closed and is called the reduction of $z$.

\begin{lemma}
For any $\omega\in\Omega_{K(X)/K}=\Omega_{K(X_\eta)/K}$, \[\RRES_{\frak{r}(z),\Cl{z}}\omega=\RRES_z\omega,\] where the left residue is the two-dimensional residue on $X$ associated to the point and curve $\frak{r}(z)\in\Cl{z}$, and the right residue is the usual residue for the $K$-curve $X_\eta$ at its closed point $z$.
\end{lemma}
\begin{proof}
This is a small exercise in chasing the definitions of the residue maps. Set $B=\roi_{X,\frak{r}(z)}$ and let $\frak{p}$ be the local equation for $\Cl{z}$ at $\frak{r}(z)$. For any $n\ge0$, $B/\frak{p}^n$ is a finite $\roi_K$-algebra, hence is complete. This implies that \[\hat{B}/\frak{p}\hat{B}=B/\frak{p},\] whence $\frak{p'}=\frak{p}\hat{B}$ is prime in $\hat{B}$, and also that \[\hat{B}_\frak{p'}/\frak{p'}^n\hat{B}_\frak{p'}=B_\frak{p}/\frak{p}^nB_\frak{p}.\] Therefore \[\hat{\hat{B}_\frak{p'}}=\projlim_n\hat{B}_\frak{p'}/\frak{p'}^n\hat{B}_\frak{p'}=\projlim_nB_\frak{p}/\frak{p}^nB_\frak{p}=\hat{B_\frak{p}}=\hat{\roi_{X_\eta,z}}.\] Then $F:=\Frac\hat{\hat{B}_\frak{p'}}$ is the two-dimensional local field used to define the residue at the flag $\frak{r}(z)\in\Cl{z}$; it has equal characteristic, and we have just shown it is equal to $\Frac\roi_{X_\eta,z}$. But the residue map on a two-dimensional local field of equal characteristic was exactly {\em defined} to be the familiar residue map for a curve.
\end{proof}

\begin{remark}\label{remark_two_dim_vs_one_dim_residues}
If $\roi_K$ is not necessarily a complete, discrete valuation ring, as at the start of the section, then the above lemma remains valid when reformulated as follows: Let $z$ be a closed point of the generic fibre, and $X_s$ a special fibre. For any $\omega\in\Omega_{K(X)/K}=\Omega_{K(X_\eta)/K}$, \[\sum_{x\in \Cl{z}\cap X_s}\RRES_{x,\Cl{z}}\omega=\RRES_z\omega\] where the left is the sum of two-dimensional residues on $X$ associated to the flags $x\in\Cl{z}$ where $x$ runs over the finitely many points in $\Cl{z}\cap X_s$, and the right residue is the usual residue at the closed point $z$ on the curve $X_\eta$. This may easily be deduced from the previous lemma using lemma \ref{lemma_residues_under_arbitrary_base_change} below.
\end{remark}

\begin{proof}[Proof of theorem \ref{theorem_reciprocity_along_a_curve}]
We may now prove the vertical reciprocity law. Let $y_1(=y),y_2,\dots,y_l$ be the irreducible components of the fibre $X_s$.

Firstly, combining the usual reciprocity law for the curve $X_\eta$ with the previous lemma yields \[\sum_{z\in (X_\eta)_0}\RRES_{\frak{r}(z),\Cl{z}}\omega=0,\] where the sum is taken over closed points of the generic fibre and only finitely many terms of the summation are non-zero. Since $\Cl{z}$, for $z\in (X_\eta)_0$, are all the irreducible horizontal curves of $X$, we may rewrite this as \[\sum_{x\in X_0}\left(\sum_{\substack{Y\subset X\sub{ horiz.}\\\sub{s.t. }Y\ni x}}\RRES_{x,Y}\omega\right)=0.\] Moreover, according to the reciprocity law around a point from subsection \ref{subsection_geometric_case}, if $x\in X_0$ is a closed point then \[\sum_{\substack{Y\subset X\\\sub{s.t. }Y\ni x}}\RRES_{x,Y}\omega=0,\] where only finitely many terms in the summation are non-zero. We deduce that \[\sum_{x\in X_0}\left(\sum_{\substack{Y\subset X\sub{ vert.}\\\sub{s.t. }Y\ni x}}\RRES_{x,Y}(\omega)\right)=0,\] where the sum is now taken over the irreducible vertical curves in $X$. That is, \[\sum_{i=1}^l\sum_{x\in y_i}\RRES_{x,y_i}\omega=0\tag{\dag},\] where the rearrangement of the double summation is justified by lemma \ref{lemma_convergence_of_residues_along_a_curve}, which says that each internal sum of (\dag) converges in $K$.

If $X_s$ is irreducible, then this is exactly the sum over the closed points of $y_1=y$ and we have finished. Else we must proceed by a `weighting' argument as in lemma \ref{lemma_continuous_wrt_m_adic_topology}. Let $\nu_1,\dots,\nu_l$ be the discrete valuations on $K(X)$ associated to $y_1,\dots,y_l$ respectively. For $m>0$, pick $f_m\in K(X)$ such that $\nu_1(f_m-1)\ge m$ and $\nu_i(f_m)\ge m$ for $i=2,\dots,l$; this exists because the $(\nu_i)_i$ are inequivalent discrete valuations. Replacing $\omega$ by $f_m\omega$ in (\dag) yields \[\sum_{i=1}^l\sum_{x\in y_i}\RRES_{x,y_i}f_m\omega=0.\] Letting $m\to\infty$ and applying the continuity part of lemma \ref{lemma_convergence_of_residues_along_a_curve} yields \[\sum_{i=1}^l\sum_{x\in y_i}\RRES_{x,y_i}f_m\omega=0\To \sum_{x\in y_1}\RRES_{x,y_1}\omega\quad\mbox{as }m\To\infty.\] This completes the proof of theorem \ref{theorem_reciprocity_along_a_curve}.
\end{proof}

\section{Trace map via residues on higher adeles}\label{section_trace_map_via_residues}
We are now ready to adelically construct Grothendieck's trace map \[H^1(X,\omegaa)\to\roi_K\] as a sum of our residues, where $\pi:X\to\Spec\roi_K$ is an arithmetic surface and $\omegaa=\omegaa_\pi$ is its relative dualizing sheaf. The key idea is to use the reciprocity laws to show that sums of residues descend to cohomology.

\begin{remark}
Passing from local constructions to global or cohomological objects is always the purpose of reciprocity laws. Compare with the reciprocity law around a point in K.~Kato and S.~Saito's two-dimensional class field theory \cite[\S4]{KatoSaito1983}. Sadly, using reciprocity laws for the reciprocity map of two-dimensional local class field theory to construct two-dimensional global class field theory has not been written down in detail anywhere, but a sketch of how it should work in the geometric case was given by Parshin \cite{Parshin1978}. More details, which are also valid in the arithmetic case, can be found in \cite[Chap.~2]{Fesenko2008a}.
\end{remark}

\subsection{Adeles of a curve}\label{subsection_adeles_of_a_curve}
We begin with a quick reminder of adeles for curves. Let $C$ be a one-dimensional, Noetherian, integral scheme with generic point $\eta$; we will be interested in both the case when $X$ is smooth over a field and when $C$ is the spectrum of the ring of integers of a number field. If $E$ is a coherent sheaf on $C$, then the {\em adelic resolution} of $E$ is the following flasque resolution:
\[0\to E\to i_\eta(E_\eta)\oplus\prod_{x\in X_0}i_x(\hat{E_x})\to\rprod_{x\in X_0}i_x(\hat{E_x}\otimes_{\roi_{X,x}} K(X))\to 0.\]
Here $i_\eta(E_\eta)$ is the constant $E_\eta$ sheaf on $X$; $\hat{E_x}$ is the $\frak{m}_{X,x}$-adic completion of $E_x$, and $i_x(\hat{E_x})$ is the corresponding skyscraper sheaf at $x$; the `restricted product' term $\rprod$ is the sheaf whose sections on an open set $U\subseteq X$ are \[\rprod_{x\in U_0}\hat{E_x}\otimes_{\roi_{X,x}}K(X)=\{(f_x)\in\prod_{x\in U_0}\hat{E_x}\otimes_{\roi_{X,x}} K(X):f_x\in\hat{E_x}\mbox{ for all but finitely many }x\in U_0\}.\] The Zariski cohomology of $E$ is therefore exactly the cohomology of the {\em adelic complex} $\bb{A}(X,E)$:
\begin{align*}
0\to E_\eta\oplus\prod_{x\in X_0}\hat{E_x}&\to\rprod_{x\in X_0}\hat{E_x}\otimes_{\roi_{X,x}}K(X)\to 0\\
(g,(f_x))&\mapsto(g-f_x)
\end{align*}
These observations remain valid if we do not bother completing $E$ at each point $x$, leading to the {\em rational adelic complex} $a(X,E)$ (classically called repartitions, e.g. \cite[II.5]{Serre1988}): \[0\to E_\eta\oplus\prod_{x\in X_0}E_x\to\rprod_{x\in X_0}E_\eta\to 0\] whose cohomology also equals the Zariski cohomology of $E$.

\subsection{Rational adelic spaces for surfaces}
The theory of adeles for curves was generalised to algebraic surfaces by A.~Parshin, e.g. \cite{Parshin1976}, and then to arbitrary Noetherian schemes by A.~Beilinson \cite{Beilinson1980}. The main source of proofs is A.~Huber's paper \cite{Huber1991}. We will describe the rational (i.e.~no completions are involved) adelic spaces, defined in \cite[\S5.2]{Huber1991}, associated to a coherent sheaf $E$ on a surface $X$. More precisely, $X$ is any two-dimensional, Noetherian, integral scheme, with generic point $\eta$ and function field $F=K(X)$. The quasi-coherent sheaf which is constantly $F$ will be denoted $\ul{F}$.

\begin{remark}
We choose to use the rational, rather than completed, adelic spaces to construct the trace map only for the sake of simplicity of notation. There is no substantial difficulty in extending the material of this section to the completed adeles, which becomes essential for the dualities discussed in remark \ref{remark_on_adelic_dualities}.
\end{remark}

\subsubsection{Adelic groups $0$, $1$, and $2$}
The first rational adelic groups are defined as follows: \[a(0)=F,\qquad a(1)=\prod_{y\in X^1}\roi_{X,y},\qquad a(2)=\prod_{x\in X^2}\roi_{X,x}.\] More generally, if $E$ is a coherent sheaf on $X$, then we define \[a(0,E)=E_\eta,\qquad a(1,E)=\prod_{y\in X^1}E_y,\qquad a(2,E)=\prod_{x\in X^2}E_x.\]

\subsubsection{Adelic group $01$}
Next we have the $01$ adelic group:
\begin{align*}
a(01)
	&=\{(f_y)\in\prod_{y\in X^1}F:\exists\mbox{ a coherent submodule }M\subseteq\ul{F}\mbox{ such that }f_y\in M_y\mbox{ for all }y\}\\
	&=\indlim_{M\subseteq\ul{F}}a(1,M)
\end{align*}
where the limit is taken over all coherent submodules $M$ of the constant sheaf $\ul{F}$. This ring is commonly denoted using restricted product notation: $a(01)=\rprod_{y\in X^1}F$. Again more generally, if $E$ is an arbitrary coherent sheaf, we put 
\begin{align*}
a(01,E)
	&=\{(f_y)\in\prod_{y\in X^1}E_\eta:\exists\mbox{ a coherent submodule }M\subseteq\ul{E_\eta}\mbox{ such that }f_y\in M_y\mbox{ for all }y\}\\
	&=\indlim_{M\subseteq\ul{E_\eta}}a(1,M),
\end{align*}
where the limit is taken over all coherent submodules $M$ of the constant sheaf associated to $E_\eta$. 

\subsubsection{Adelic group $02$}
Next,
\begin{align*}
a(02)
	&=\{(f_x)\in\prod_{x\in X^2} F:\exists\mbox{ a coherent submodule }M\subseteq\ul{F}\mbox{ such that }f_x\in M_x\mbox{ for all }x\}\\
	&=\indlim_{M\subseteq\ul{F}}a(2,M),
\end{align*}
where the limit is taken over all coherent submodules $M$ of $\ul{F}$. This ring is commonly denoted $\rprod_{x\in X^2}F$. We leave it to the reader to write down the definition of $a(02,E)$, for $E$ an arbitrary coherent sheaf.

\subsubsection{Adelic group $12$}
\begin{remark}
We first require some notation. If $z\in X$ is any point and $N$ is a $\roi_{X,z}$ module, then we write \[[N]_z=j_{z*}(\tilde{N}),\] where $j_z:\Spec\roi_{X,z}\into X$ is the natural morphism and $\tilde{N}$ is the quasi-coherent sheaf on $\Spec\roi_{X,z}$ induced by $N$. For example, $\ul{F}=[\roi_{X,\eta}]_\eta$.
\end{remark}

We may now introduce
\[a(12)=\prod_{y\in X^1}a_y(12),\] where 
\begin{align*}
a_y(12)
	&=\{(f_x)\in\prod_{x\in y} \roi_{X,y}:\exists\mbox{ a coherent submodule }M\subseteq[\roi_{X,y}]_y\mbox{ such that }f_x\in M_x\mbox{ for all }x\in y\}\\
	&=\indlim_{M\subseteq[\roi_{X,y}]_y}a(2,M),
\end{align*}
where the limit is taken over all coherent submodules $M$ of $[\roi_{X,y}]_y$. Recall our convention that if $y\in X^1$ then `$x\in y$' means that $x$ is a codimension one point of the closure of $y$; more precisely, $x\in X^2\cap\res{\{y\}}$.

We again leave it to the reader to write down the definition of $a(12,E)$ for an arbitrary coherent sheaf $E$ (just replace $\roi_{X,y}$ by $E_y$ everywhere in the construction).

This is a convenient place to make one observation concerning an adelic condition which holds for $a(12,E)$:

\begin{lemma}\label{lemma_12_adelic_condition2}
Let $E$ be a coherent sheaf on $X$, fix $y\in X^1$, $r\ge 0$, and let $(f_x)_{x\in y}\in a_y(12, E)$; then $f_x\in E_x+\frak{m}_{X,y}^rE_y$ for all but finitely many $x\in y$.
\end{lemma}
\begin{proof}
There is a coherent submodule $M\subseteq[E_y]_y$ such that $f_x\in M_x$ for all $x\in y$. Let $U=\Spec A$ be an affine open neighbourhood of $y$, and let $\frak{q}\subset A$ be the prime ideal corresponding to $y$. Then $M(U)$ is a finitely generated $A$-submodule of $E_\frak{q}$ and therefore $M(U)\subseteq fE$ for some $f\in A_\frak{q}$. For any $r\ge 0$, the argument of lemma \ref{lemma_12_adelic_condition} shows that $f\in A_\frak{m}+\frak{q}^rA_\frak{q}$ for all but finitely many of the maximal ideals $\frak{m}$ of $A$ sitting over $\frak{q}$; for such maximal ideals, we have $M_\frak{m}\subseteq E_\frak{m}+\frak{q}^rE_\frak{q}$. Since $U$ contains all but finitely many of the points of $\Cl y$, this is enough.
\end{proof}

\subsubsection{Adelic group $012$}
Finally, \[a(012)=\indlim_{M\subseteq\ul{F}}a(12,M)\subseteq\prod_{y\in X^1}\prod_{x\in y}F.\] (and we similarly define $a(012, E)$ for any coherent $E$, by taking the limit over coherent submodules $M$ of the constant sheaf $\ul{E_\eta}$).

\subsubsection{Simplicial structure and cohomology}
Consider the following homomorphisms of rings:
\[\xymatrixcolsep{1cm}\xymatrixrowsep{1cm}\xymatrix{
& 	& 	F\ar[dl]\ar[dr]&	&	\\
&	\prod_{y\in X^1}F\ar[r] & \prod_{y\in X^1}\prod_{x\in y}F & \prod_{x\in X^2}F\ar[l] & \\
\prod_{y\in X^1}\roi_{X,y}\ar[rr]\ar[ur] & & \prod_{y\in X^1}\prod_{x\in y}\roi_{X,y}\ar[u] & & \prod_{x\in X^2}\roi_{X,x}\ar[ll]\ar[ul]
}\]
where the three ascending arrows are the obvious inclusions and the remaining arrows are diagonal embeddings. These homomorphisms restrict to the rational adelic groups just defined to give a commutative diagram of ring homomorphisms:
\[\xymatrixcolsep{1.5cm}\xymatrixrowsep{1.5cm}\xymatrix{
& 	& 	a(0)\ar[dl]_{\bor^0_{01}}\ar[dr]^{\bor^0_{02}}&	&	\\
&	a(01)\ar[r]^{\bor^{01}_{012}} & a(012) & a(02)\ar[l]_{\bor^{02}_{012}} & \\
a(1)\ar[rr]_{\bor^1_{12}}\ar[ur]^{\bor^1_{01}} & & a(12)\ar[u]_{\bor^{12}_{012}} & & a(2)\ar[ll]^{\bor^2_{12}}\ar[ul]_{\bor^2_{02}}
}\] (and similarly with any coherent sheaf $E$ in place of $\roi_X$). For example, to see that $\bor^1_{12}$ is well-defined, once must check that if $f\in\roi_{X,y}$ then there is a coherent submodule $M$ of $[\roi_{X,y}]_y$ such that $f_x\in M_x$ for all $x\in y$; but $f$ may be viewed as a global section of $[\roi_{X,y}]_y$ and therefore $M:=\roi_Xf\subseteq[\roi_{X,y}]_y$ suffices.

We reach the analogue for $X$ of the rational adelic complex which we saw for a curve in subsection \ref{subsection_adeles_of_a_curve} above:

\begin{theorem}\label{theorem_adeles_compute_cohomology}
Let $E$ be a coherent sheaf on $X$; then the Zariski cohomology of $E$ is equal to the cohomology of the following complex (which is the total complex associated to the above simplicial group):
\comment{
\[\xymatrixcolsep{0mm}\xymatrixrowsep{0.4cm}\xymatrix{0\ar[r]&
a(0,E)\oplus a(1,E)\oplus a(2,E)\ar[r]& a(01,E)\oplus a(02,E)\oplus a(12,E)\ar[r]& a(012,E)\ar[r]&0\\
&(f_0,f_1,f_2)\ar@{|->}[r]& (\bor^{0}_{01}f_0-\bor^{1}_{01}f_1, \bor^{2}_{02}f_2-\bor^{0}_{02}f_0,\bor^{1}_{12}f_1-\bor^{2}_{12}f_2) &&\\
&&(g_{01},g_{02},g_{12})\ar@{|->}[r]&\bor^{01}_{012}g_{01}+\bor^{02}_{012}g_{02}+\bor^{12}_{012}g_{12}&}\]
}
\begin{align*}
&0\To a(0,E)\oplus a(1,E)\oplus a(2,E)\To a(01,E)\oplus a(02,E)\oplus a(12,E)\To a(012,E)\To0\\ \\
&\phantom{0\To a(0,E)}
(f_0,f_1,f_2)\mapsto (\bor^{0}_{01}f_0-\bor^{1}_{01}f_1, \bor^{2}_{02}f_2-\bor^{0}_{02}f_0,\bor^{1}_{12}f_1-\bor^{2}_{12}f_2)\\\\
&\phantom{0\To a(0,E)\oplus a(1,E)\oplus a(2,E)\To a(01,E)}
(g_{01},g_{02},g_{12})\mapsto\bor^{01}_{012}g_{01}+\bor^{02}_{012}g_{02}+\bor^{12}_{012}g_{12}
\end{align*}

\end{theorem}
\begin{proof}
This is due to Parshin \cite{Parshin1976}; the general case of higher dimensional $X$ is due to Beilinson \cite{Beilinson1980} and Huber \cite{Huber1991}.
\end{proof}

\subsection{Construction of the trace map}
Let $\roi_K$ be a Dedekind domain of characteristic zero with finite residue fields; its field of fractions is $K$. Let $\pi:X\to S=\Spec\roi_K$ be an $\roi_K$-curve as at the start of section \ref{section_reciprocity_along_vertical_curves}. According to the main result of \cite{Morrow2009}, the relative dualising sheaf $\omegaa$ of $\pi$ is explicitly given by, for open $U\subseteq X$, \begin{align*}\omegaa(U)=\{\omega\in\Omega_{K(X)/K}:&\RRES_{x,y}(f\omega)\in\hat{\roi_{K,\pi(x)}}\mbox{ for}\tag{\dag}\\& \mbox{all } x\in y\subset U\mbox{ and }f\in\roi_{X,y}\}\end{align*} where $x$ runs over all closed points of $X$ inside $U$ and $y$ runs over all curves containing $x$.

As previously, closed points of $S$ are denoted $s$, and we put $\roi_s=\hat{\roi_{K,s}}$, $K_s=\Frac\roi_s$.

\begin{proposition}\label{proposition_residues_respect_adelic_conditions}
If $\ul{\omega}=(\omega_{x,y})_{x\in y}\in a(012,\omegaa)$ and $s\in S_0$, then \[\RRES_s(\ul{\omega}):=\sum_{\substack{x,y\\\sub{s.t. }x\in y\cap X_s}}\RRES_{x,y}\omega_{x,y}\tag{\ddag}\] converges in $K_s$, where the sum is taken over all points $x$ and curves $y$ in $X$ for which $x\in y\cap X_s$. Moreover, $\RRES_s(\ul{\omega})\in\roi_s$ for all but finitely many $s\in S_0$.

If $\ul{\omega}\in\bor^{12}_{012}a(12,\omegaa)$ then all terms of the sum, hence also $\RRES_s(\ul{\omega})$, belong to $\roi_s$.
\end{proposition}
\begin{proof}
Let $E$ be a coherent submodule of the constant sheaf $\ul{\omegaa}_\eta=\ul{\Omega}_{K(X)/K}$ such that $\ul{\omega}\in a(12,E)$; then $E$ and $\omegaa$ are equal at the generic point (replacing $E$ by $E+\omegaa$, if necessary), hence on an open set, and therefore $E_y=\omegaa_y$ for all but finitely many $y\in X^1$. We call the remaining finitely many $y$ {\em bad}.

If $y$ is a horizontal curve which is not bad and $x\in y$, then $\omega_{x,y}\in E_y=\omegaa_y$ and so $\RRES_{x,y}\omega_{x,y}=0$ (indeed, if $\pi\in\roi_{K,s}$ is a uniformiser at $s$ then $\pi^{-1}\in\roi_{X,y}$ and so the definition of $\omegaa$ implies that $\pi^{-m}\RRES_{x,y}\omega_{x,y}\in\roi_s$ for all $m\ge0$; this is only possible if $\RRES_{x,y}\omega_{x,y}=0$). Therefore, only finitely many horizontal curves contribute to the summation in (\ddag); so it is enough to prove that if $y$ is an irreducible component of $X_s$ then \[\sum_{x\in y} \RRES_{x,y}\omega_{x,y}\] converges. This is straightforward, using lemma \ref{lemma_12_adelic_condition2} and arguing exactly as in lemma \ref{lemma_convergence_of_residues_along_a_curve}, and completes the proof that $\RRES_s(\ul{\omega})$ is well-defined.

Secondly, for any curve $y$, each of $\omegaa_y$ and $E_y$ are (non-zero) finitely generated $\roi_{X,y}$ submodules of $\Omega_{K(X)/K}$, and therefore there exists $r\ge0$ such that $\frak{m}_{X,y}^rE_y\subseteq\omegaa_y$; clearly we may pick $r$ so that this inclusion holds for all bad $y$. Then lemma \ref{lemma_12_adelic_condition2} tells us that for all but finitely many $x$ in any bad curve $y$, we have \[E_y\subseteq E_x+\frak{m}_{X,y}^rE_y\subseteq E_x+\omegaa_y.\] Next, if $y_1, y_2$ are two horizontal curves, then $y_1$ and $y_2$ will have a common point of intersection on a vertical curve $Y$ for only finitely many $Y$ (for else $y_1\cap y_2$ would be infinite). It follows that there is an open set $U\subseteq X$ consisting of fibres such that any $x\in U$ satisfies one of the following conditions:
\begin{enumerate}
\item $x$ sits on no bad curve; or
\item $x$ sits on exactly one bad curve $y$; $y$ is horizontal and $E_y\subseteq E_x+\omegaa_y$.
\end{enumerate}
Note that $U$ contains all but finitely many of the fibres $X_s$, for $s\in S_0$, and to prove our second claim it is enough to show that for any closed point $x$ on a fibre $X_s$ belonging to $U$, and curve $y$ passing through $x$, one has $\RRES_{x,y}\omega_{x,y}\in\roi_s$. There are two cases to consider:
\begin{enumerate}
\item $y$ is not bad. Then $\omega_{x,y}\in E_y=\omegaa_y$, whence $\RRES_{x,y}\omega_{x,y}\in\roi_s$ by (\dag).
\item $y$ is bad. Then $y$ is horizontal by construction of $U$ and so $\RRES_{x,y}\omegaa_y=0$ (as argued in the previous paragraph); therefore condition (ii) on $U$ implies that $\RRES_{x,y}\omega_{x,y}=\RRES_{x,y}\zeta$ for some $\zeta\in E_x$. If $Y$ is any curve through $x$ apart from $y$ then $\zeta\in E_x\subseteq E_Y=\omegaa_Y$ and so (\dag) now implies that $\RRES_{x,Y}\zeta\in\roi_s$. But the reciprocity law about a point from subsection \ref{subsection_geometric_case} shows that \[\RRES_{x,y}\zeta=-\sum_Y\RRES_{x,Y}\zeta,\] where the sum is taken over all curves $Y$ passing through $x$ apart from $y$; therefore $\RRES_{x,y}\zeta\in\roi_s$.
\end{enumerate}
This completes the proof that $\RRES_s\ul{\omega}$ belongs to $\roi_s$ for all but finitely many $s\in S_0$.

Finally, if $\ul{\omega}$ is in the image of the boundary map $\bor^{12}_{012}$ then $\omega_{x,y}\in\omegaa_y$ for all flags $x\in y$; so (\dag) implies that $\RRES_{x,y}\omega_{x,y}\in\roi_s$. This proves the final claim.
\end{proof}

Let \[\bb{A}_S=\rprod_{s\in S_0}K_s=\{(a_s)\in\prod_{s\in S_0}K_s:a_s\in\roi_s\mbox{ for all but finitely many }s\}\] and \[\bb{A}_S(0)=\prod_{s\in S_0}\roi_s\] be the rings of adeles and integral adeles of $K$ respectively (we will incorporate archimedean information in the final section). The adelic complex for $S$, as discussed in \ref{subsection_adeles_of_a_curve}, is
\begin{align*}
0\To K\oplus\bb{A}_S(0)&\To\bb{A}_S\To 0\\
(\lambda,(a_s))&\mapsto (\lambda-a_s)
\end{align*}

\begin{corollary}
The map \[\RRES:a(012,\omegaa)\to\bb{A}_S,\quad \ul{\omega}\mapsto(\RRES_s(\ul{\omega}))_{s\in S_0}\] is well-defined, and restricts to $\RRES\circ\bor^{12}_{012}:a(12,\omegaa)\to\bb{A}_S(0)$.
\end{corollary}
\begin{proof}
This is exactly the content of the previous proposition.
\end{proof}

Define a map
\begin{align*}
\RRES':a(01,\omegaa)\oplus a(02,\omegaa)\oplus a(12,\omegaa&)\to K\oplus\bb{A}_S(0)\\
(\ul{\omega}',\ul{\omega}'',\ul{\omega})\phantom{a(12,\omega)}&\mapsto (\sum_{z\in X_\eta}\RRES_z\omega_z', \RRES(\bor^{12}_{012}\ul{\omega}))
\end{align*}
where the first sum is taken over closed points $z$ of $X_\eta$ or, equivalently, horizontal curves in $X$, and $\RRES_z$ denotes the usual residue for $X_\eta$ as a smooth curve over $K$ (note that this makes sense as $\omegaa_\eta=\Omega_{K(X_\eta)/K}$). In the remainder of the paper, $z$ will always denote a closed point of $X_\eta$.

The key application of the reciprocity laws is to deduce that taking sums of residues induces a morphism of adelic complexes:

\begin{proposition}
The following maps give a homomorphism of adelic complexes from $X$ to $S$:
\[\xymatrixcolsep{0.5cm}\xymatrixrowsep{1cm}\xymatrix{
0\ar[r]& a(0,\omegaa)\oplus a(1,\omegaa)\oplus a(2,\omegaa)\ar[r]\ar[d]& a(01,\omegaa)\oplus a(02,\omegaa)\oplus a(12,\omegaa)\ar[r]\ar[d]^{\RRES'}& a(012,\omegaa)\ar[r]\ar[d]^{\RRES} &0\\
& 0 \ar[r] & K\oplus\bb{A}_S(0) \ar[r] & \bb{A}_S \ar[r] & 0
}\]
\end{proposition}
\begin{proof}
Commutativity of the first square is equivalent to the following results:
\begin{enumerate}
\item If $\omega\in a(0,\omegaa)=\Omega_{K(S)/K}$ then $\sum_{z\in X_\eta}\RRES_z\omega=0$.
\item If $\ul{\omega}=(\omega_y)_{y\in X^1}\in a(1,\omegaa)$ then $\sum_{z\in X_\eta}\RRES_z\omega_z=0$ and $\RRES(\bor^{12}_{012}\bor^1_{12}\ul{\omega})=0$.
\item If $\ul{\omega}\in a(2,\omegaa)$ then $\RRES(\bor^{12}_{012}\bor^2_{12}\ul{\omega})=0$.
\end{enumerate}

(i) is the usual reciprocity law for the curve $X_\eta/K$. The first vanishing claim in (ii) holds since $\omega_z\in\omegaa_z=\Omega_{X_\eta/K,z}$ and the residue of a differential form on $X_\eta$ at a point where it is regular is zero. For the second vanishing claim in (ii), note that if $s\in S_0$ then \[\RRES_s(\bor^{01}_{012}\bor^1_{01}\ul{\omega})=\sum_{y\subseteq X_s}\sum_{x\in y}\RRES_{x,y}\omega_y + \sum_{\substack{\sub{horiz.}\\y}}\sum_{x\in X_s\cap y}\RRES_{x,y}\omega_y,\] where we have split the summation (\ddag) depending on whether $y$ is an irreducible component of $X_s$ or is horizontal. But the first double summation is zero, according to the reciprocity law along a vertical curve (theorem \ref{theorem_reciprocity_along_a_curve}), while every term in the second double summation is zero since they are residues along horizontal curves $y$ of forms in $\omegaa_y$ (see the second paragraph of the previous proof). We will return to (iii) in a moment.

Commutativity of the second square is almost automatic since $\RRES'$ was obtained by restricting $\RRES$ to $a(01,\omegaa)$ and $a(12,\omegaa)$; it remains only to check that if $\ul{\omega}\in a(02,\omegaa)$ then $\RRES\bor^{02}_{012}\ul{\omega}=0$. This follows immediately from the reciprocity law around a point from \ref{subsection_geometric_case}. This also establishes (iii), since if $\ul{\omega}\in a(2,\omegaa)$ then $\bor^{12}_{012}\bor^2_{12}\ul{\omega}=\bor^{02}_{012}\bor^2_{02}\ul{\omega}\in\bor^{02}_{012}a(02,\omegaa)$.
\end{proof}

Noting that $H^0$ of the adelic complex for $S$ is simply $\roi_K$ and that $H^1$ of the adelic complex for $X$ is $H^1(X,\omegaa)$ (by theorem \ref{theorem_adeles_compute_cohomology}), the proposition implies that there is an induced map \[\RRES:H^1(X,\omegaa)\to \roi_K.\] Our construction would be irrelevant without the final theorem:

\begin{theorem}
$\RRES$ is equal to Grothendieck's trace map $\operatorname{tr}_\pi$.
\end{theorem}
\begin{proof}
There is a natural morphism from the rational adelic complex of $X$ for the coherent sheaf $\omegaa$ to the rational adelic complex of $X_\eta$ for the coherent sheaf $\Omega_{X_\eta/K}$:
\[\xymatrixcolsep{0.8cm}\xymatrixrowsep{1cm}\xymatrix{
0\ar[r]& a(0,\omegaa)\oplus a(1,\omegaa)\oplus a(2,\omegaa)\ar[r]\ar[d]^{(\omega_0,\omega_1,\omega_2)\mapsto (\omega_0,p_1(\omega_1))}& a(01,\omegaa)\oplus a(02,\omegaa)\oplus a(12,\omegaa)\ar[r]\ar[d]^{(\omega_{01},\omega_{02},\omega_{12})\mapsto p_{01}(\omega_{01})}& a(012,\omegaa)\ar[r]\ar[d] &0\\
0 \ar[r] & \Omega_{K(X)/K}\oplus\prod_{z\in X_\eta}\Omega_{X_\eta/K,z} \ar[r] & \rprod_{z\in X_\eta}\Omega_{K(X)/K} \ar[r] & 0 \\
}\]
This is given by the identity $a(0,\omegaa)=\Omega_{K(X)/K}$, the projection
\[a(1,\omegaa)=\prod_{y\in X^1}\omegaa_y=\prod_{z\in X_\eta}\Omega_{X_\eta/K,z}\times\prod_{\substack{y\in X^1\\\sub{vertical}}}\omegaa_y\stackrel{p_1}{\onto}\prod_{z\in X_\eta}\Omega_{X_\eta/K,z},\]
and the restriction of the projection
\[\prod_{y\in X^1}\omegaa_\eta=\prod_{z\in X_\eta}\Omega_{K(X)/K}\times\prod_{\substack{y\in X^1\\\sub{vertical}}}\Omega_{K(X)/K}\onto \prod_{z\in X_\eta}\Omega_{K(X)/K}\]
to the adelic spaces
$a(01,\omegaa)\stackrel{p_2}{\onto}\rprod_{z\in X_\eta}\Omega_{K(X)/K}$.

By the functoriality of adeles, the resulting map $H^*(X,\omegaa)\to H^*(X_\eta,\Omega_{X_\eta/K})$ is the natural map on cohomology induced by the restriction $\omegaa|_{X_\eta}=\Omega_{X_\eta/K}$. Using this, we will now show that
\[\xymatrix{
H^1(X,\omegaa) \ar[r]\ar[d]_{\RRES} & H^1(X_\eta,\Omega_{X_\eta/K})\ar[d]^{\operatorname{tr}}\\
\roi_K\ar[r] & K 
}\tag{\maltese}\]
commutes, where the right vertical arrow is the trace map for the $K$-curve $X_\eta$. Indeed, from the definition of $\RRES'$ above, the following diagram certainly commutes:
\[\xymatrixcolsep{4cm}\xymatrix{
\ker\langle a(01,\omegaa)\oplus a(02,\omegaa)\oplus a(12,\omegaa)\to a(012,\omegaa)\rangle \ar[d]_{\RRES'}\ar[r]^-{(\omega_{01},\omega_{02},\omega_{12})\mapsto p_{01}(\omega_{01})} & \rprod_{z\in X_\eta}\Omega_{K(X)/K}\ar[d]_{(\omega_z)\mapsto\sum_{z\in X_\eta}\RRES_z\omega_z}\\
\ker\langle K\oplus\bb{A}_S(0)\to\bb{A}_S\rangle=\roi_K\ar[r] & K
}\]
Passing to cohomology groups, we deduce that 
\comment{\[\xymatrix{
H^1(X,\omegaa) \ar[r]\ar[d]_{\RRES} & H^1(X_\eta,\Omega_{X_\eta/K})=\operatorname{Coker}\langle \Omega_{K(X)/K}\oplus\prod_{z\in X_\eta}\Omega_{X_\eta/K,z} \to \rprod_{z\in X_\eta}\Omega_{K(X)/K}\rangle \ar[d]^{(\omega_z)\mapsto\sum_{z\in X_\eta}\RRES_z\omega_z}\\
\roi_K\ar[r] & K
}\]}
\[\xymatrix{
H^1(X,\omegaa) \ar[r]\ar[d]_{\RRES} & H^1(X_\eta,\Omega_{X_\eta/K}) \ar[d]^{(\omega_z)\mapsto\sum_{z\in X_\eta}\RRES_z\omega_z} \ar@{}[r]|-{=}& \operatorname{Coker}\langle \Omega_{K(X)/K}\oplus\prod_{z\in X_\eta}\Omega_{X_\eta/K,z} \to \rprod_{z\in X_\eta}\Omega_{K(X)/K}\rangle \\
\roi_K\ar[r] & K & 
}\]
commutes; but the vertical map on the right is the trace map for $X_\eta$, by the familiar result (which we are generalising!) that the trace map of a smooth projective curve is represented by the sum of residues. This completes the proof that (\maltese) commutes.

Finally, the diagram (\maltese) also commutes if $\RRES$ is replaced by $\operatorname{tr}_\pi$, since trace maps commute with localization of the base ring. Therefore $\RRES=\operatorname{tr}_\pi$.
\end{proof}

\begin{remark}\label{remark_higher_dimension}
Before complicating matters by incorporating archimedean data, this is a convenient opportunity to explain how the previous material should fit into a general framework.

A {\em flag of points} on a scheme $X$ is a sequence of points $\xi=(x_0,\dots,x_n)$ such that $x_{i-1}\in\res{\{x_i\}}$ for $i=1,\dots,n$. By a process of successive completions and localizations, the flag $\xi$ yields a ring $F_\xi$. More generally, to any quasi-coherent sheaf $E$, one obtains a module $E_\xi$ over $F_\xi$; for details, see \cite[\S3.2]{Huber1991}.

Now let $f:X\to Y$ be a morphism of $S$-schemes, where $S$ is a Noetherian scheme (perhaps Cohen-Macaulay), and notice that we may push forward any flag from $X$ to $Y$, \[f_*(\xi):=(f(x_0),\dots,f(x_n)),\] resulting in an inclusion of rings $F_{f_*(\xi)}\subseteq F_\xi$. Let $\omegaa_X$, $\omegaa_Y$ denote the dualizing sheaves of $X$, $Y$ over $S$. If $f$ is proper (and probably Cohen-Macaulay) of fibre dimension $d$, then we expect there to exist a residue map \[\RRES_\xi:\omegaa_{X,\xi}\to\omegaa_{Y,f_*(\xi)}\] which is the trace map when $f$ is a finite morphism and which is transitive when given another proper, CM morphism $Y\to Z$. Globally, taking sums of these residue maps will induce a morphism of degree $-d$ on the adelic complexes \[\RRES_{X/Y}:\bb{A}(X,\omegaa_X)\to\bb{A}(Y,\omegaa_Y).\] The patching together of the local residue maps to induce a morphism of complexes is equivalent to a collection of reciprocity laws being satisfied. In turn, this induces maps on the cohomology \[H^*(X,\omegaa_X)=H^*(\bb{A}(X,\omegaa_X))\To H^{*-d}(\bb{A}(Y,\omegaa_Y))=H^{*-d}(Y,\omegaa_Y),\] which will be nothing other than Grothendieck's trace map.

When $S$ is a field this framework more-or-less follows from \cite{Lomadze1981} and \cite{Yekutieli1992}, though it has not been written down carefully. This article and the author's previous \cite{Morrow2009} focus on the case $Y=S=\Spec\roi_K$ and $X=$ a surface.

The fully general case requires a rather careful development of relative residue maps in arbitrary dimensions, and becomes a technically difficult exercise quite quickly. The Hochschild homology-theoretic description of residue maps \cite{Hubl1989} \cite{Lipman1987} may be the key to a smoother approach.
\end{remark}

\section{Archimedean reciprocity along horizontal curves}\label{section_archimedean}
We continue to study an $\roi_K$-curve $X$ in the sense introduced at the start of section \ref{section_reciprocity_along_vertical_curves}, but we now assume that $K$ is a number field and $\roi_K$ its ring of integers (with generic point $\eta$). If $\infty$ is an infinite place of $K$ then we write $X_\infty=X\times_{\roi_K} K_\infty$ where $K_\infty$ is the completion of $K$ at $\infty$; so $X_\infty$ is a smooth projective curve over $\bb{R}$ or $\bb{C}$.

The natural morphism \[X_\infty=X\times_{\roi_K} K_\infty\xto{\rho} X_\eta= X\times_{\roi_K}K\] can send a closed point to the generic point; but there are only finitely many points over any closed point. Indeed, let $z\in X_\eta$ be a closed point; then the fibre over $z$ is \[X_\infty\times_{X_\eta} k(z)=(K_\infty\times_K X_\eta)\times_{X_\eta} k(z)=\Spec (K_\infty\otimes_K k(z)),\] which is a finite reduced scheme.

If $y$ is a horizontal curve on $X$ then $y=\Cl{z}$ for a unique closed point $z\in X_\eta$. We say that a closed point $x\in X_\infty$ sits on $y$ if and only if $\rho(x)=z$. Hence there are only finitely many points on $X_\infty$ which sit on $y$, and we will allow ourselves to denote this set of points by $y\cap X_\infty$. Such points are the primes of $K_\infty\otimes_K k(z)$ and therefore correspond to the infinite places of the number field $k(z)$ extending the place $\infty$ on $K$. Note that each $x\in X_\infty$ sits on at most one horizontal curve, which may seem strange at first.

In this situation, we define the {\em archimedean residue map} $\RRES_{x,y}:\Omega_{K(X)/K}\to K_\infty$ to be \[\Omega_{K(X)/K}\To\Omega_{K(X_\infty)/K_\infty}\xto{\RRES_x}K_\infty,\] where $\RRES_x$ is the usual one-dimensional residue map associated to the closed point $x$ on the smooth curve $X_\infty$ over $K_\infty$.

The following easy lemma was used in remark \ref{remark_two_dim_vs_one_dim_residues}; since we need it again, let's state it accurately:

\begin{lemma}\label{lemma_residues_under_arbitrary_base_change}
Let $C$ be a smooth, geometrically connected curve over a field $K$ of characteristic zero, let $L$ be an arbitrary extension of $K$, and let $z$ be a closed point of $C$.
\begin{enumerate}
 \item Let $x\in C_L$ be a closed point sitting over $z$; then the following diagram commutes:
\[\xymatrixcolsep{2cm}\xymatrix{
\Omega_{K(C_L)/L}\ar[r]^{\RES_x}&k(x)\\
\Omega_{K(C)/K}\ar[r]^{\RES_z}\ar[u]&k(z)\ar[u]\\
}\]
(Notation: $\RES_x$ is the residue map $\Omega_{K(C)/K}\to k(x)$, and then $\RRES_x=\Tr_{k(x)/K}\circ\RES_x$; similarly for other points.)
\item With $x$ now varying over all the closed points of $C_L$ sitting over $z$, the following diagram commutes:
\[\xymatrixcolsep{2cm}\xymatrix{
\Omega_{K(C_L)/L}\ar[r]^{\sum_{x|z}\RRES_x}&L\\
\Omega_{K(C)/K}\ar[r]^{\RRES_z}\ar[u]&K\ar[u]\\
}\]
\end{enumerate}
\end{lemma}
\begin{proof}
If $t\in K(C)$ is a local parameter at $z$ then it is also a local parameter at $x$, and the characteristic zero assumption implies that there are compatibile isomorphisms $K(C_L)_x\cong k(x)((t))$, $K(C)_z\cong k(z)((t))$; the first claim easily follows. Secondly $k(z)\otimes_KL\cong\bigoplus_{x|z}k(x)$, so that $\Tr_{k(z)/K}=\sum_{x|z}\Tr_{k(x)/L}$; hence, for $\omega\in\Omega_{K(C)/K}$, part (i) lets us use the usual argument:
\begin{align*}
\sum_{x|z}\RRES_x(\omega)
	&=\sum_{x|z}\Tr_{k(x)/L}\RES_x(\omega)\\
	&=\sum_{x|z}\Tr_{k(x)/L}\RES_z(\omega)\\
	&=\Tr_{k(z)/K}\RES_z(\omega)\\
	&=\RRES_z(\omega)\qedhere
\end{align*}
\end{proof}

We obtain an analogue of remark \ref{remark_two_dim_vs_one_dim_residues}: 

\begin{corollary}
Returning to the notation before the lemma, if $\infty$ and $y=\Cl{z}$ are fixed, and $\omega\in\Omega_{K(X)/K}$, then \[\sum_{x\in y\cap X_\infty}\RRES_{x,y}\omega=\RRES_z\omega.\]
\end{corollary}
\begin{proof}
Apply the previous lemma with $C=X_\eta$ and $L=K_\infty$.
\end{proof}

Write $\res{S}=\Spec\roi_K\cup\{\infty\mbox{'s}\}$ for the `compactification' of $S=\Spec\roi_K$ by the infinite places (in fact, the notation $s\in\res{S}$ will always mean that $s$ is a place of $K$, never the generic point of $S$) and let \[\bb{A}_{\res{S}}=\rprod_{s\in\res{S}}K_s=\bb{A}_S\times\prod_{\infty}K_\infty\] be the usual ring of adeles of the number field $K$. Let \[\psi=\otimes_{s\in\res{S}}\psi_s:\bb{A}_{\res{S}}\to S^1\quad\mbox{($=$ the circle group\footnotemark)}\] \footnotetext{{We never consider the set of codimension one points of $S=\Spec\roi_K$, so this shouldn't cause confusion.}} be a continuous additive character which is trivial on the global elements $K\subset\bb{A}_{\res{S}}$ \cite[Lem.~4.1.5]{Tate1967}.

Note that, if $y$ is a horizontal curve on $X$, then even with our definition of points at infinity, it does not make sense consider a reciprocity law \[\mbox{``}\sum_{x\in y}\RRES_{x,y}\omega=0\mbox{''}\] since the residues appearing live in different local fields. This problem is fixed by using the ``absolute base'' $S^1$:

\begin{definition}
Let $y$ be a curve on $X$ and $x\in y$ a closed point sitting over $s\in\res{S}$ (this includes the possibility that $y$ is horizontal and $s$ is an infinite place). Define the {\em absolute residue map} \[\psi_{x,y}:\Omega_{K(X)/K}\to S^1\] to be the composition \[\Omega_{K(X)/K}\xto{\RRES_{x,y}}K_s\xto{\psi_s}S^1.\]
\end{definition}

We may now establish the reciprocity law on $X$ along {\em any} curve, including the horizontal ones:

\begin{theorem}
Let $y$ be a curve on $X$ and $\omega\in\Omega_{K(X)/K}$. Then for all but finitely many closed points $x\in y$ the absolute residue $\psi_{x,y}(\omega)$ is $1$, and \[\prod_{x\in y}\psi_{x,y}(\omega)=1\] in $S^1$.
\end{theorem}
\begin{proof}
First consider the case that $y$ is an irreducible component of a special fibre $X_s$ (here $s\in S_0$). Then $\ker\psi_s$ is an open subgroup of $K_s$, and so the proof of lemma \ref{lemma_convergence_of_residues_along_a_curve} shows that $\RRES_{x,y}\omega\in\ker\psi_s$ for all but finitely many $x\in y$. Also, \[\prod_{x\in y}\psi_{x,y}(\omega)=\psi_s\left(\sum_{x\in y}\RRES_{x,y}(\omega)\right),\] which is $\psi_s(0)=1$ according to the reciprocity law along the vertical curve $y$ (theorem \ref{theorem_reciprocity_along_a_curve}).

Secondly suppose that $y=\Cl{z}$ is a horizontal curve; here $z$ is a closed point of $X_\eta$. The proof of proposition \ref{proposition_residues_respect_adelic_conditions} shows that $\RRES_{x,y}\omega\in\roi_{\pi(x)}$ for all but finitely many $x\in y$ (here $x$ is a genuine schematic point on $X$); since $\ker\psi_s$ contains $\roi_s$ for all but finitely many $s\in S_0$, it follows that $\psi_{x,y}(\omega)=1$ for all but finitely many $x\in y$. It also follows that \[\ul{f}:=\left(\sum_{x\in y\cap X_s}\RRES_{x,y}\omega\right)_{s\in\res{S}}\] belongs to $\bb{A}_{\res{S}}$, and clearly \[\prod_{x\in y}\psi_{x,y}(\omega)=\prod_{s\in\res{S}}\psi_s\left(\sum_{x\in y\cap X_s}\RRES_{x,y}\omega\right)=\psi(\ul{f}).\] But remark \ref{remark_two_dim_vs_one_dim_residues} (for $s\in S_0$) and the previous corollary (for $s$ infinite) imply that $\ul{f}$ is the global adele $\RRES_z\omega\in K$. As $\psi$ was chosen to be trivial on global elements, the proof is complete.
\end{proof}

\begin{remark}
The reciprocity law around a point $x\in X^2$ stated in section \ref{subsection_geometric_case} obviously implies that the absolute residue maps satisfy a similar law: \[\prod_{\substack{y\subset X\\\sub{s.t. }y\ni x}}\psi_{x,y}(\omega)=1.\]

Therefore we have absolute reciprocity laws for all points and for all curves, which are analogues for an arithmetic surface of the reciprocity laws established by Parshin \cite{Parshin1976} for an algebraic surface.
\end{remark}

\begin{remark}\label{remark_on_adelic_dualities}
Let $F_{x,y}$ be the finite direct sum of two-dimensional local fields attached to a flag $x\in y$; i.e.~$F_{x,y}=\Frac\hat{A_\frak{p}}$ where $A=\hat{\roi_{X,x}}$, $\frak{p}=y\roi_{X,x}$, and $y\subset\roi_{X,x}$ also denotes the local equation for $y$ at $x$; so $F_{x,y}=\bigoplus_{y'|y} F_{x,y}$ where $y'$ varies over the finitely many height one primes of $A$ over $y$, and $F_{x,y'}=\Frac\hat{A_{y'}}$.

By the local construction of the residue maps we see that $\psi_{x,y}$ is really the composition \[\Omega_{K(X)/K}\To\Omega_{K(X)/K}\otimes_{K(X)}F_{x,y}=\bigoplus_{y'|y}\Omega_{F_{x,y'}/K_s}^\sub{cts}\xto{\sum_{y'|y}\RRES_{F_{x,y'}}}K_s\xto{\psi_s}S^1,\] and each $\psi_s\circ\RRES_{F_{x,y'}}:F_{x,y'}\to S^1$ is a continuous (with respect to the two-dimensional topology; see remark \ref{remark_continuous_wrt_two_dim_topology}) character on the two-dimensional local field $F_{x,y'}$. This character will induce self-duality of the topological group $F_{x,y'}$, which in turn will induce various dualities on the (complete) adelic groups; for some results in this direction, see \cite[\S27, \S28]{Fesenko2008a}. 
\end{remark}

\begin{remark}
Taking $S=\Spec\bb{Z}$, it would be very satisfying to have an extension of the framework discussed in remark \ref{remark_higher_dimension} to include archimedean points. The main existing problem is the lack at present of a good enough theory of adeles in arbitrary dimensions which include the points at infinity. The author is currently trying to develop such a theory and hopes that this will allow the dualities discussed in the previous remark to be stated more precisely and in greater generality (in all dimensions and including points at infinity).
\end{remark}

\bibliographystyle{acm}
\bibliography{../Bibliography}

\noindent Matthew Morrow,\\
University of Chicago,\\
5734 S. University Ave.,\\
Chicago,\\
IL, 60637,\\
USA\\
{\tt mmorrow@math.uchicago.edu}\\
\url{http://math.uchicago.edu/~mmorrow/}\\
\end{document}